\documentclass{amsart}

\usepackage[mathcal]{euscript}
\usepackage{amssymb, amsfonts, xypic}
\usepackage{enumerate}

\newtheorem{theorem}{Theorem}[section]
\newtheorem{lemma}[theorem]{Lemma}
\newtheorem{corollary}[theorem]{Corollary}
\newtheorem{proposition}[theorem]{Proposition}

\theoremstyle{definition}
\newtheorem{definition}[theorem]{Definition}
\newtheorem{remark}[theorem]{Remark}
\newtheorem{example}[theorem]{Example}

\numberwithin{equation}{section}

\newcommand{\Z}{\mathbb{Z}}
\newcommand{\N}{\mathbb{N}}
\newcommand{\Cx}{\mathbb{C}}
\newcommand{\R}{\mathbb{R}}

\newcommand{\C}{\mathcal{C}}

\newcommand{\D}{\mathcal{D}}

\DeclareMathOperator{\Hom}{Hom}

\DeclareMathOperator*{\colim}{colim}
\DeclareMathOperator*{\holim}{holim}

\newcommand{\pros}{\textup{pro}}
\newcommand{\pro}{\textup{pro-}}
\newcommand{\Map}{\textup{Map}}

\newcommand{\map}{\rightarrow}

\newcommand{\tl}{\tilde}

\newcommand{\bdry}[1]{\partial \Delta[#1]}

\newcommand{\dfn}{\textbf} %Make defined words bold.
\newcommand{\mdfn}[1]{\dfn{\mathversion{bold}#1}} % Even make math symbols bold

\begin{document}

\title{Generalized cohomology of pro-spectra}

\author{Daniel C. Isaksen}
\address{Department of Mathematics \\ Wayne State University\\
Detroit, MI 48202}

\email{isaksen@math.wayne.edu}

\date{July 30, 2004}

\thanks{The author was partially supported 
by an NSF Postdoctoral Fellowship.  The
author thanks Tom Nevins, Halvard Fausk, Steve Costenoble, and Stefan Waner
for suggestions and encouragement.}

\begin{abstract}
We present a closed model structure for the category of pro-spectra
in which the weak equivalences are detected by stable homotopy pro-groups.
With some bounded-below assumptions, weak equivalences are also detected
by cohomology as in the classical Whitehead theorem for spectra.
We establish an Atiyah-Hirzebruch spectral sequence in this context,
which makes possible the computation of topological $K$-theory 
(and other generalized cohomology theories) of pro-spectra.
\end{abstract}

\maketitle

\section{Introduction}
\label{sctn:intro}

Ordinary singular cohomology of pro-objects is a useful tool in various
mathematical concepts.  For example, the cohomology of pro-spaces comes
into play in the Bousfield-Kan viewpoint on $R$-completions of spaces
\cite{BK} \cite{D} \cite{I4}.  Also, the singular cohomology 
with locally constant coefficients of the 
\'etale homotopy type of a scheme \cite{AM} \cite{F2}
is isomorphic to the \'etale cohomology
of the scheme.
The continuous cohomology of a pro-finite group \cite{S} is also
an example of the same kind.

The notion of the 
cohomology of a pro-object is easy to describe.
For any cofiltered system $X$, we define
\begin{equation}
\label{eq:cohomology}
H^{*} (X) = \colim_{s} H^{*} (X_{s}).
\end{equation}

Wherever singular cohomology is useful, it is a good bet that generalized
cohomology theories are also useful.  This paper develops the foundations 
and tools necessary for studying these cohomology theories on pro-objects.
In particular, we give a definition of generalized cohomology
theories such that there is an Atiyah-Hirzebruch spectral
sequence whose convergence is reasonably well-behaved 
(see Theorem \ref{thm:AHSS}).

We are only aware of 
one example of generalized cohomology theories applied to
pro-objects: \'etale $K$-theory,
which is the topological $K$-theory of
the \'etale homotopy type of a scheme \cite{F1}.  
In future work, we plan to
use the foundations in this paper to develop a homotopy fixed points spectral
sequence for pro-finite groups \cite{FI}.  We also intend to use
\'etale $BP$ and \'etale $tmf$, which are analogous to \'etale $K$-theory,
in order to prove some results concerning quadratic forms over fields
of characteristic $p$ \cite{DI1} \cite{DI2}.

Unfortunately, generalized cohomology theories are not as easy to define
as ordinary cohomology.
In view of Formula (\ref{eq:cohomology}),
the most obvious definition of the topological $K$-theory
of a pro-object is the formula
\[
KU^*(X) = \colim_{s} KU^{*} (X_{s}).
\]
However, this turns out to be wrong computationally.
Since there is an Atiyah-Hirzebruch spectral sequence whose abutment is
each $KU^*(X_s)$ and since filtered colimits are exact, 
one might hope that the filtered colimit of these
spectral sequences gives a computational tool for understanding $KU^*(X)$
as defined above.
Unfortunately, the convergence of this colimit spectral sequence is terrible
and therefore of no practical use.  See Section \ref{sctn:K-theory} 
for a specific
computational example of the problem.

The correct way to define generalized cohomology theories is more complicated.
First, we define a closed model structure on the 
category of pro-spectra (see Section \ref{sctn:model-structure}).  
This gives a homotopy theory
for pro-spectra and in particular a homotopy category.
Then, for any spectrum $E$,
we simply define $E^r(X)$ to be the set of homotopy
classes of maps from $\Sigma^{-r} X$ to the spectrum $E$ considered as a 
constant pro-spectrum.
In the same way, we can define a cohomology theory represented by
any pro-spectrum.

Roughly speaking,
the weak equivalences of the model structure on pro-spectra are 
defined in terms of pro-homotopy groups 
(see Section \ref{sctn:pro-homotopy}),
and the fibrant pro-spectra are cofiltered diagrams of spectra whose
homotopy groups are bounded above (see Section \ref{sctn:fibrant}).
This characterization of the fibrant pro-spectra
connects with the fact 
that the Postnikov tower of the $K$-theory spectrum $KU$ 
is key to the definition of \'etale $K$-theory \cite{DF}. 
This paper gives a concrete explanation 
for why the Postnikov tower enters the picture because
the Postnikov tower of $KU$, thought of as a pro-spectrum, is a fibrant
replacement for $KU$.

The model structure for pro-spectra is analogous
to a model structure for pro-spaces \cite{I1}.
Many of the technical complications of \cite{I1}
are questions of choosing basepoints and studying $\pi_{1}$-actions.
These issues do not arise for pro-spectra, so many aspects of the
theory for pro-spectra are easier than for pro-spaces.

There are at least two other reasonable ways to put a model structure on the
category of pro-spectra \cite{EH} \cite{CI}.  Each gives a distinct
homotopy category, and each has particular uses.  
We have two strong pieces of evidence that the model structure under
investigation in this paper is the correct one for the study of 
generalized cohomology theories of pro-spectra.

First, there is a Whitehead theorem showing that weak equivalences
of pro-spectra can be detected by cohomology isomorphisms under
some bounded below hypotheses (see Theorem \ref{thm:Whitehead}).
This is relevant
in many applications because often pro-objects are constructed precisely
for their cohomological properties.

The second piece of evidence is the Atiyah-Hirzebruch spectral sequence
discussed above whose convergence is well-behaved.

\subsection{Organization}

The paper is organized as follows.  
We start with a concrete example in Section \ref{sctn:K-theory}
to demonstrate what is wrong with the naive definition of topological
$K$-theory of a pro-spectrum.

Section \ref{sctn:prelim-pro} is a review of the
machinery of pro-categories and the very general ``strict'' model
structure \cite[\S~3.3]{EH} \cite{I2}.  
The strict model structure is the starting
point for all known homotopy theories of pro-categories. 

Then we recall in 
Section \ref{sctn:prelim-spectra}
some ideas about spectra and stable homotopy theory.
All of the basic properties of spectra that we will need are 
satisfied by all of the usual models of spectra, such as
Bousfield-Friedlander spectra \cite{BF},
symmetric spectra \cite{HSS}, or $S$-modules \cite{EKMM}.
Thus, the results in this paper can be viewed as applying to any of these
categories of spectra.

Then in Section \ref{sctn:model-structure} 
we define the cofibrations, fibrations, and
weak equivalences of pro-spectra, and we
prove that they are a model structure.
We assume that the reader has a basic familiarity with the terminology
and standard results of model categories.
The original source is \cite{Q}, but 
we follow the notation and terminology of \cite{Hi}
as closely as possible.  Other references include 
\cite{Ho} and \cite{DS}.

The next few sections contain useful properties of the homotopy
theory of pro-spectra.  
If one is ever going to use the model structure on pro-spectra for anything,
it is essential to know what the cofibrant and fibrant objects are.
The cofibrant objects are easy to describe.
In Section \ref{sctn:fibrant}, we identify explicitly the fibrant
pro-spectra.  

In Section \ref{sctn:homotopy-class} we 
collect some results about computing homotopy classes
of maps of pro-spectra in terms of homotopy classes of maps of spectra.

The actual definition of weak equivalences (see Definition
\ref{defn:weak-equivalence}) has the advantage that it is useful
for proving model structure axioms.  However, it lacks a computational
aspect.
In Section \ref{sctn:pro-homotopy}, we make precise the relationship
between pro-homotopy groups and weak equivalences of pro-spectra.
This tends to be useful in applications.
In Section \ref{sctn:Whitehead}, 
we give yet another description of the weak equivalences 
in terms of cohomology.

Finally, in Section \ref{sctn:AHSS} we construct the Atiyah-Hirzebruch
spectral sequence and prove that it is conditionally convergent
in the sense of \cite{B} for a large class of pro-spectra 
of interest.  For example, the pro-spectra $\R P^{\infty}_{-\infty}$
and $\Cx P^{\infty}_{-\infty}$ \cite{L} both belong to this class.

There are two ways to construct the Atiyah-Hirzebruch spectral sequence
in ordinary stable homotopy theory converging to $[X,Y]^{*}$.
One uses the skeletal filtration of $X$, and the other uses the 
Postnikov tower of $Y$ \cite[App.~B]{GM}.  We use the approach
with Postnikov towers here.
This is no surprise because the Postnikov towers play such
an important role in the model structure for pro-spectra.

We do not address the question of multiplicative properties of
the Atiyah-Hirze\-bruch spectral sequence for pro-spectra, but we strongly
suspect that everything works as expected.

\section{What is the $K$-theory of a pro-spectrum?}
\label{sctn:K-theory}

We take the position that however the $K$-theory of a pro-spectrum is
defined, there ought to be an Atiyah-Hirzebruch spectral sequence
with reasonably good convergence properties.  If $K$-theory is to be
useful computationally, this is an appropriate expectation.

Let us assume for the moment that the $K$-theory of a pro-spectrum
$X = \{ X_s \}$ is defined to be $\colim_s KU^*(X_s)$.  We will show that this
definition does lead to an Atiyah-Hirzebruch spectral sequence, but
the convergence is not at all good.

For each $s$, there is an Atiyah-Hirzebruch spectral sequence
\[
H^{p}(X_s; \pi_{q} KU) \Rightarrow KU^{p+q}(X_s).
\]
Since filtered colimits are exact, we can take colimits and obtain
a spectral sequence
\[
\colim_s H^p(X_s; \pi_q KU) \Rightarrow \colim_s KU^{p+q}(X_s).
\]
The left side is just the ordinary cohomology
$H^p(X; \pi_q KU)$ of $X$ with coefficients in $\pi_q KU$, 
so this appears to be an Atiyah-Hirzebruch spectral sequence.

It remains to ask what kind of convergence properties this spectral
sequence has.  In the following particular case, we will show that
the convergence is terrible; it's so bad that the $E_2$-term basically
gives no information at all about the abutment.

For each $n \geq 0$, let $X_n$ be the spectrum 
$\displaystyle\vee_{k=n}^\infty S^{2k}$,
and let $X_{n+1} \map X_n$ be the obvious inclusion.  Thus $X$ is
a pro-spectrum (in fact, a countable tower).

Since $X_n$ is $(2n-1)$-connected, $H^p(X_n; \Z)$ is zero for sufficiently
large $n$.  Therefore, $\colim_n H^p(X_n;\pi_q KU)$ is zero for all
$p$ and $q$.  This means
that the $E_2$-term of the above spectral sequence is zero.

On the other hand, $KU^{p+q}(X_n)$ is equal to 
$\prod_{k=n}^\infty \Z$ when $p+q$ is even and
equal to zero when $p+q$ is odd.  Therefore, when $p+q$ is even,
$\colim_n KU^{p+q}(X_n)$
is equal to a quotient of $\prod_{k=1}^\infty \Z$, where
two infinite sequences $(a_k)$ and $(b_k)$ are identified if 
$a_k$ and $b_k$ are different for only finitely many values of $k$.
Another way to think of this group is the ``germs at infinity'' of
functions $\N \map \Z$.

When $p+q$ is even, $\colim_n KU^{p+q}(X_n)$ is uncountable.  Recall
that the $E_2$-term of the spectral sequence was zero.  The conclusion
is that the above spectral sequence has disastrously bad convergence
properties.

We are led to the conclusion that $\colim_s KU^*(X_s)$ is the wrong
definition of the $K$-theory of a pro-spectrum.  The point of the 
rest of this paper is to construct a suitable homotopy theory of
pro-spectra such that $[X, KU]_{\pros}$ does have the desired computational
properties.  Here, $KU$ means the constant pro-spectrum with value $KU$.
When we define the weak equivalences in this homotopy theory later,
it will be clear that for $X$ in the previous paragraphs,
the map $* \map X$ is a weak equivalence.  Therefore,
$[X,KU]_{\pros}$ is necessarily zero, which agrees with the computation
of the $E_2$-term of the spectral sequence.

\section{Preliminaries on Pro-Categories} 
\label{sctn:prelim-pro}

We begin with a brief review of pro-categories.
This section contains mostly standard material on pro-categories
\cite{SGA}
\cite{AM} \cite{EH}.
We conform to the notation and terminology of \cite{I2}.

\subsection{Pro-Categories}
\label{subsctn:pro}

\begin{definition} 
\label{defn:pro}
For a category $\C$, the category \mdfn{$\pro \C$} has objects all 
cofiltering
diagrams in $\C$, and 
$$\Hom_{\pro \C}(X,Y) = \lim_s \colim_t \Hom_{\C}
     (X_t, Y_s).$$
Composition is defined in the natural way.
\end{definition}

A \dfn{constant} pro-object is one indexed
by the category with one object and one (identity) map.
Let $\mdfn{c}: \C \map \pro \C$ be the functor taking an object $X$ to the 
constant pro-object with value $X$.
Note that this functor makes $\C$ a full subcategory of
$\pro \C$.  The limit functor $\lim: \pro \C \map \C$
is the right adjoint of $c$.  

\subsection{Level representations}
\label{subsctn:level}

A \dfn{level map} $X \map Y$ is a pro-map that is given
by a natural transformation (so $X$ and $Y$ must have the same indexing
category); this is very special kind of pro-map.  
A \dfn{level representation} of a pro-map
$f:X \map Y$ is
another pro-map $\tl{f}: \tl{X} \map \tl{Y}$ such that 
$\tl{f}$ is a level map.
Moreover, we require that there are isomorphisms
$X \map \tl{X}$ and $Y \map \tl{Y}$ such that the diagram
\[
\xymatrix{
X \ar[d] \ar[r]^{f} & Y \ar[d] \\
\tl{X} \ar[r]_{\tl{f}} & \tl{Y}     }
\]
of pro-maps commutes.
Every map has a level representation 
\cite[App.~3.2]{AM}.  See \cite[App.]{C} for a functorial construction
of level representations.

A pro-object $X$ satisfies a certain property \dfn{levelwise} if
each $X_{s}$ satisfies that property, and $X$ satisfies this property
\dfn{essentially levelwise} if it is isomorphic to another pro-object
satisfying this property levelwise.
Similarly, a level map $X \map Y$ 
satisfies a certain property \dfn{levelwise}
if each $X_s \map Y_s$ has this property.
A map of pro-objects satisfies this property \dfn{essentially levelwise}
if it has a
level representation satisfying this property levelwise.

The following purely technical lemma will be needed later 
in Lemma \ref{lem:detect-acyclic-cofibrations}.

\begin{lemma}
\label{lem:pro-isomorphism}
Let $Y$ be a pro-object.
Suppose that for some of the maps $t \map s$ in the indexing diagram
for $Y$, there exists an object $Z_{ts}$ and a factorization
$Y_t \map Z_{ts} \map Y_s$ of the structure map $Y_t \map Y_s$.
Also suppose that for every $s$, there exists at least one $t \map s$
with this property.
The objects $Z_{ts}$ assemble into a pro-object $Z$ that is isomorphic
to $Y$.
\end{lemma}

\begin{proof}
We may assume that $Y$ is indexed by a directed set $I$ (in the sense
that there is at most one map between any two objects of $I$)
because every pro-object is isomorphic to a pro-object indexed by
a directed set \cite[Thm.~2.1.6]{EH}.
Define a new directed set $K$ as follows.  The elements of 
$K$ consist of pairs $(t,s)$ of elements of $I$ such that
$t \geq s$ and a factorization $Y_t \map Z_{ts} \map Y_s$ exists.
If $(t,s)$ and $(t',s')$ are two elements of $K$, we say that
$(t',s') \geq (t,s)$ if $s' \geq t$.  It can easily
be checked that this makes $K$ into a directed set.

Note that the function $K \map I: (t,s) \mapsto s$
is cofinal in the sense of \cite[App.~1]{AM}.  
This means that we may reindex
$Y$ along this functor and assume that $Y$ is indexed by $K$;
thus we write $Y_{(t,s)} = Y_s$.

We define the pro-spectrum $Z$ to be indexed by $K$ by setting
$Z_{(t,s)} = Z_{ts}$.  If $(t',s') \geq (t,s)$, then the structure map
$Z_{(t',s')} \map Z_{(t,s)}$ is the composition
\[
Z_{t's'} \map Y_{s'} \map Y_{t} \map Z_{ts},
\]
It can easily be checked that this gives a functor defined on $K$;
here is where we use that the composition $Y_t \map Z_{ts} \map Y_s$
equals $Y_t \map Y_s$.

Finally, we must show that $Z$ is isomorphic to $Y$.  We use the 
criterion from \cite[Lem.~2.3]{I1} for detecting pro-isomorphisms.
Given any $(t,s)$ in $K$, choose $u$ such that $(u,t)$ is in $K$.
Then there exists a diagram
\[
\xymatrix{
Z_{(u,t)} \ar[d]\ar[r] & Y_{(u,t)} = Y_t \ar[d]\ar[dl] \\
Z_{(t,s)} \ar[r] & Y_{(t,s)} = Y_s.                     }
\]
\end{proof}

\subsection{Strict Model Structures} 
\label{subsctn:strict}

We make some remarks on the \dfn{strict model structure} for pro-categories,
originally developed in \cite{EH} and studied further in \cite{I2}.
The niceness hypothesis of \cite[\S~2.3]{EH} is not satisfied by
the categories of spectra that we will use, so the generalizations of
\cite{I2} really are necessary.
The categories of pro-simplicial sets, pro-topological spaces, and any of the
standard models for pro-spectra (such as Bousfield-Friedlander spectra
\cite{BF}, symmetric spectra \cite{HSS}, or $S$-modules \cite{EKMM}) 
all have strict model structures.

Let $\C$ be a proper model category.
The \dfn{strict weak equivalences} 
of $\pro \C$ are the essentially
levelwise weak equivalences (see Section \ref{subsctn:level}).
The \dfn{cofibrations} of $\pro \C$ are the
essentially levelwise cofibrations.  Finally, the 
\dfn{strict fibrations} of
$\pro \C$ are maps that have the right lifting property with respect
to the strict acyclic cofibrations.
We use no adjective to describe the cofibrations because the cofibrations
are the same in all known model structures on pro-categories.

The following theorem is the main result of \cite{I2}.

\begin{theorem} \label{thm:strict}
Let $\C$ be a proper model category.  Then the classes of
cofibrations, strict weak equivalences, and strict fibrations
define a proper model structure on $\pro \C$.  If $\C$ is simplicial,
then this structure is also simplicial.
\end{theorem}

For any two objects $X$ and $Y$ of $\pro \C$,
the mapping space $\Map(X,Y)$ is equal to $\lim_s \colim_t \Map(X_t, Y_s)$ 
when $\C$ is simplicial.

We will need the following fact in a few places.  It makes computations
of mapping spaces significantly easier.

\begin{proposition} \label{prop:strict-map}
Let $\C$ be a proper simplicial model category.
Let $X$ be a cofibrant object of $\pro \C$.  
Let $Y$ be any levelwise fibrant object of $\pro \C$
with strict fibrant replacement $\hat{Y}$.
Then the homotopically
correct mapping space $\Map(X,\hat{Y})$ is weakly equivalent to
$\holim_s \colim_t \Map(X_t, Y_s)$.
\end{proposition}

\begin{proof}
Since $\Map(X, \hat{Y})$ is homotopically correct, it doesn't matter
which strict fibrant replacement $\hat{Y}$ that we consider.  Therefore,
we may choose one with particularly good properties.  Use the method
of \cite[Lem.~4.7]{I2} to factor the map $Y \map *$ into a strict
acyclic cofibration $Y \map \hat{Y}$ followed by a strict fibration
$\hat{Y} \map *$.  This particular construction gives that 
$Y \map \hat{Y}$ is a levelwise weak equivalence and that $\hat{Y}$
is levelwise fibrant.

Since $X$ is cofibrant and $\hat{Y}$ is strict fibrant,
the pro-space $s \mapsto \colim_t \Map(X_t, \hat{Y}_s)$ 
is also strict fibrant.
This can be seen by inspecting the explicit description of 
strict fibrations given in \cite[Defn.~4.2]{I2}.
Therefore, $\Map(X, \hat{Y}) = \lim_s \colim_t \Map(X_t, \hat{Y}_s)$ 
is weakly equivalent to
$\holim_s \colim_t \Map(X_t, \hat{Y}_s)$ because
homotopy limit is the derived functor of limit with respect
to the strict model structure \cite[Rem.~4.2.11]{EH}.

The map $\colim_t \Map(X_t, Y_s) \map \colim_t \Map(X_t, \hat{Y}_s)$ 
is a weak equivalence because
$Y_s \map \hat{Y}_s$ is a weak equivalence between fibrant objects.
Homotopy limits preserve levelwise weak equivalences, so the map
\[
\holim_s \colim_t \Map(X_t, Y_s) \map \holim_s \colim_t \Map(X_t, \hat{Y}_s)
\]
is a weak equivalence.
\end{proof}

We will next show that construction of the strict model structure
respects Quillen equivalences \cite[Defn.~8.5.20]{Hi}.  
It was an oversight that
this result was not included in \cite{I2}.

If $F: \C \map \D$ is any functor, then there is another functor
$F: \pro \C \map \pro \D$ defined by applying $F$ levelwise to
any object in $\pro \C$.  If $G: \D \map \C$ is right adjoint to $F$,
then $G$ is also right adjoint to $F$ on pro-categories.

\begin{theorem}
\label{thm:strict-Quillen-equivalence}
Let $\C$ and $\D$ be model categories such that the strict model
structures on $\pro \C$ and $\pro \D$ exist (for example, if
$\C$ and $\D$ are proper).
If $F: \C \map \D$ and $G: \D \map \C$ are a Quillen adjoint pair,
then $F$ and $G$ are also
a Quillen adjoint pair between $\pro \C$ and $\pro \D$ equipped
with their strict model structures.  If $F$ and $G$ are a Quillen
equivalence between $\C$ and $\D$, then they are also a Quillen
equivalence between $\pro \C$ and $\pro \D$.
\end{theorem}

\begin{proof}
First suppose that $F$ and $G$ are a Quillen adjoint pair on $\C$ and $\D$.
Since $F$ takes cofibrations in $\C$ to cofibrations in $\D$, it
takes levelwise cofibrations in $\pro \C$ to levelwise cofibrations
in $\pro \D$.  Thus, $F$ preserves
essentially levelwise cofibrations.

Since $F$ takes acyclic cofibrations in $\C$ to
acyclic cofibrations in $\D$, it similarly preserves essentially
levelwise acyclic cofibrations.  However, the essentially levelwise
acyclic cofibrations are the same as the strict acyclic 
cofibrations \cite[Prop.~4.11]{I2}.
This shows that $F$ preserves strict acyclic cofibrations. 
Thus $F$ and $G$ are a Quillen adjoint pair.

Now suppose that $F$ and $G$ are a Quillen equivalence on $\C$ and $\D$.
To show that $F$ and $G$ are a Quillen equivalence on $\pro \C$
and $\pro \D$,
let $X$ be a cofibrant object of $\pro \C$ and let $Y$ be
a strict fibrant object of $\pro \D$.
Suppose that $g:X \map GY$ is a strict weak equivalence; we want to show
that its adjoint 
$f: FX \map Y$ is also a strict weak equivalence.

We may assume that $X$ is levelwise cofibrant.  
By \cite[Lem.~4.5]{I2}, we may also assume
that $Y$ is levelwise fibrant.  
Using the level replacement of
\cite[App.~3.2]{AM}, we may reindex $X$ and $Y$ in such a
way that $X$ is still levelwise cofibrant, $Y$ is still levelwise fibrant,
and $g: X \map GY$ is a level map.
However, we are not allowed to assume that $g$ is
a levelwise weak equivalence because this may require a different
reindexing.

Use the method of \cite[Lem.~4.6]{I2} to factor $f$ into
a strict cofibration $i: X \map Z$ followed by a
strict acyclic fibration $p:Z \map GY$.  This particular construction gives
that $i$ is a levelwise cofibration and that $p$ is a levelwise
acyclic fibration.  In particular, this implies that $Z$ is levelwise cofibrant
since $X$ is.  The two-out-of-three axiom implies that
$i$ is a strict acyclic cofibration, even though it is not a levelwise
acyclic cofibration.

The adjoint $p': FZ \map Y$ of $p$ is a levelwise weak equivalence
because $F$ and $G$ are a Quillen equivalence between $\C$ and $\D$.
This works because $Z$ is levelwise cofibrant, $Y$ is levelwise fibrant,
and $p$ is a level map.
The map $Fi: FX \map FZ$ is a strict acyclic cofibration because left Quillen
functors preserve acyclic cofibrations.
The map $f$
is the composition of $Fi$ with $p'$, so $f$ is a strict weak equivalence.

Now assume that $f: FX \map Y$ is a strict weak equivalence.  To show
that its adjoint $g: X \map GY$ is also a strict weak equivalence,
use the dual argument.
\end{proof}

In particular, $- \wedge S^1$ and $\Map(S^1, -)$ are adjoint functors
on spectra.  Thus, they induce adjoint functors on pro-spectra also.
Similarly to spectra, we define the suspension functor $\Sigma$ 
and loops functor $\Omega$ on pro-spectra
as the derived functors of these functors.
We won't recall the basic details
of spectra until the next section.  For now, it is enough to know that
the functors $- \wedge S^1$ and $\Map(S^1, -)$ are a Quillen equivalence
from the category of spectra to itself.

\begin{theorem}
\label{thm:strict-stable}
The functors $- \wedge S^1$ and $\Map(S^1, -)$ are a Quillen equivalence
from the strict model structure on pro-spectra to itself.
\end{theorem}

In other words, the strict model structure on pro-spectra
is stable in the sense of \cite{Ho}.  

\begin{proof}
This is an immediate application of 
Theorem \ref{thm:strict-Quillen-equivalence}.
\end{proof}

\section{Preliminaries on Spectra} 
\label{sctn:prelim-spectra}

This section contains some results on spectra and stable homotopy theory.
Much of the material is well-known.

We work with a proper simplicial model structure on a category of spectra
such as Bousfield-Friedlander spectra \cite{BF}, $S$-modules \cite{EKMM},
or symmetric spectra \cite{HSS}.
We assume that the model structure is cofibrantly generated.
Moreover, the cofiber of any generating cofibration must be a sphere.
If the dimension of the sphere is $k$, then
we call such a map a \mdfn{generating cofibration of dimension $k$}.  

We also need that stable weak equivalences are detected by stable homotopy
groups (this is true even for symmetric spectra if the stable homotopy
groups are properly defined).  Also, stable homotopy groups commute with
colimits along transfinite compositions of cofibrations.

As usual, 
the symbol \mdfn{$\Sigma$} refers to the \dfn{suspension} functor,
the left derived version of the functor $- \wedge S^1$.
Thus $\Sigma X$ is defined to be $\tilde{X} \wedge S^1$ for a cofibrant 
replacement $\tilde{X}$ of $X$.
Similarly, the symbol \mdfn{$\Omega$} refers to the \dfn{loops} functor,
the right derived version of the functor $\Map(S^1, -)$.
This means that $\Omega X$ is defined to be $\Map(S^1,\hat{X})$ 
for a fibrant replacement $\hat{X}$ of $X$.
The key property of $\Sigma$ and $\Omega$ is that they
are inverse equivalences on the stable homotopy category.

Let $[X, Y]$ be the set of stable homotopy classes from $X$ to $Y$, and let
$[X, Y]^{r}$ be the set of stable 
homotopy classes of degree $r$ from $X$ to $Y$.  
If the functor $\Sigma^{r}$ is defined to be $\Omega^{-r}$ for
$r \leq 0$, then
$[X, Y]^{r}$ is equal to 
\[
[\Sigma^{-r} X, Y] = [X, \Sigma^{r} Y]
\]
for all $r$.

An \dfn{Eilenberg-Mac Lane spectrum} is a spectrum whose stable homotopy
groups are zero except in one dimension.

\subsection{$n$-Equivalences}
\label{subsctn:n-equivalence}

In the next two subsections, we study some special kinds of maps of 
spectra that play a central role in the model structure for pro-spectra.

\begin{definition} \label{defn:n-equivalence}
A map $f$ of spectra is an \mdfn{$n$-equivalence} if $\pi_{k} f$ is an
isomorphism for $k < n$ and $\pi_{n} f$ is a surjection.  A map
$f$ is a \mdfn{co-$n$-equivalence} if $\pi_{k} f$ is an
isomorphism for $k > n$ and $\pi_{n} f$ is an injection.
\end{definition}

\begin{definition} \label{defn:bounded}
A spectrum $X$ is \dfn{bounded below} if the map $* \map X$ is an
$n$-equivalence for some $n$.
A spectrum $X$ is \dfn{bounded above} if the map $X \map *$ is a
co-$n$-equivalence for some $n$.
\end{definition}

Of course, a bounded below spectrum is a spectrum whose homotopy
groups vanish below some (arbitrarily small) dimension,
and a bounded above spectrum is a spectrum whose homotopy
groups vanish above some (arbitrarily large) dimension.

\begin{lemma} \label{lem:n-equivalence}
A map
is an $n$-equivalence if and only if its homotopy cofiber $C$
satisfies $\pi_{k} C = 0$ for all $k \leq n$.
A map 
is a co-$n$-equivalence if and only if its homotopy fiber $F$
satisfies $\pi_{k} F = 0$ for all $k \geq n$.
\end{lemma}

\begin{proof}
This follows immediately from the long exact sequence of homotopy
groups of a homotopy cofiber sequence or homotopy fiber sequence.
\end{proof}

\begin{lemma} \label{lem:pullback-pushout-n-equivalence}
Base changes along fibrations preserve 
$n$-equivalences and co-$n$-equi\-va\-lences.
Cobase changes along cofibrations preserve $n$-equivalences
and co-$n$-equi\-va\-lences.
\end{lemma}

\begin{proof}
First consider a pullback square
\[
\xymatrix{
W \ar[r] \ar[d]_g & Z \ar[d]^{f} \\
X \ar[r]_{p} & Y}
\]
in which $p$ is a fibration and $f$ is a co-$n$-equivalence.
Let $F$ be the homotopy fiber of $f$.  By Lemma \ref{lem:n-equivalence},
$\pi_k F = 0$ for all $k \geq n$.  The pullback square is a homotopy pullback
square
because $p$ is a fibration, so the homotopy fiber of $g$ is also $F$.
By Lemma \ref{lem:n-equivalence} again, $g$ is a co-$n$-equivalence.

Now suppose that $f$ is an $n$-equivalence.
Then the homotopy cofiber $C$ of $f$ satisfies $\pi_k C = 0$ for all
$k \leq n$.  Note that $C$ is the suspension $\Sigma F$ of the 
homotopy fiber $F$.
Since $F$ is also the homotopy fiber of $g$, $C$ is also the
homotopy cofiber of $g$.
Lemma \ref{lem:n-equivalence} again tells us that $g$ is an
$n$-eqvuivalence.

The proof for cobase changes along cofibrations is dual.
\end{proof}

\subsection{Co-$n$-Fibrations and $n$-Cofibrations}
\label{subsctn:n-cofibration-fibration}

Now we need some results on how the $n$-equivalences interact with the
fibrations and cofibrations.

\begin{definition}\label{defn:n-cofibration-fibration}
A map of spectra is 
a \mdfn{co-$n$-fibration} 
if it has the right lifting property with respect
to all generating acyclic cofibrations and
all generating cofibrations of dimension greater than $n$.
A map of spectra is
an \mdfn{$n$-cofibration} if it has the 
left lifting property with respect to all co-$n$-fibrations.
\end{definition}

Note that co-$n$-fibrations and
$n$-cofibrations are characterized by lifting properties with
respect to
each other.
Also, the class of $n$-cofibrations is the same as the class of
retracts of relative $J_{n}$-cell complexes, 
where \mdfn{$J_{n}$} is the set of 
generating acyclic cofibrations together with the set of generating
cofibrations of dimension greater than $n$
\cite[Cor.~10.5.23, Defn.~12.4.7]{Hi}.

When $n = -\infty$,
the definitions reduce to the usual definitions of cofibrations
and acyclic fibrations.
When $n = \infty$,
the definitions reduce to the usual definitions of acyclic cofibrations
and fibrations.

\begin{lemma} \label{lem:n-cofibration-fibration}
Every acyclic fibration is a co-$n$-fibration, and
every co-$n$-fibration is a fibration.  
Every acyclic cofibration is an $n$-cofibration, and every 
$n$-cofibration is a cofibration.
If $m \geq n$, then every $m$-cofibration is an $n$-cofibration,
and every co-$n$-fibration is a co-$m$-fibration.
\end{lemma}

\begin{proof}
Compare the lifting properties given in Definition 
\ref{defn:n-cofibration-fibration} to the usual lifting properties
of cofibrations, acyclic cofibrations, fibrations, and acyclic fibrations.
\end{proof}

\begin{lemma}\label{lem:n-factor}
For any $n$, maps of spectra factor functorially
into $n$-cofibrations followed by co-$n$-fibrations.
\end{lemma}

\begin{proof}
Apply the small object argument \cite[Prop.~10.5.16]{Hi} 
to the set $J_{n}$ of acyclic generating cofibrations together with
generating cofibrations of dimension greater than $n$.
\end{proof}

When working with $n$-cofibrations and co-$n$-fibrations,
we use the following two propositions frequently to pass between
lifting properties and properties of homotopy groups as expressed
in the notions of $n$-equivalences and co-$n$-equivalences.

\begin{proposition} \label{prop:co-n-fibration}
A map of spectra is a co-$n$-fibration if and only if it is a fibration
and a co-$n$-equivalence.
\end{proposition}

\begin{proof}
This is proved in \cite[Thm.~8.6]{CDI}.  Here is the basic idea.
Obstructions for lifting generating cofibrations
of dimension $k$ with respect to a fibration $p$
belong to the $(k-1)$st stable homotopy group of the fiber of $p$.
This connects to co-$n$-equivalences via 
Lemma \ref{lem:n-equivalence}.
\end{proof}

\begin{proposition} \label{prop:n-cofibration}
A map $f$ of spectra is an $n$-cofibration if and only if it is a cofibration
and an $n$-equivalence.
\end{proposition}

\begin{proof}
Consider the class $C$ of all maps that are cofibrations and $n$-equivalences.
We will first show that $C$ contains all retracts of $J_n$-cell complexes and
thus contains all $n$-cofibrations.

Acyclic cofibrations belong to $C$, as do generating cofibrations
of dimension greater than $n$.  Therefore $C$ contains $J_n$.  
An argument similar to the proof of 
Lemma \ref{lem:pullback-pushout-n-equivalence} implies that $C$
is closed under cobase changes.  Next, observe that $C$ is closed
under transfinite compositions of cofibrations 
because stable homotopy groups commute with filtered colimits along such 
compositions.  Finally, 
retracts preserve cofibrations
and $n$-equivalences, so $C$ is closed under retracts.
This finishes one implication.  

For the other implication, assume that
$f$ is a cofibration and $n$-equivalence.  Let $C$ be the cofiber of $f$.
By Lemma \ref{lem:n-equivalence},
the desuspension $\Omega C$ has the property that 
$\pi_k \Omega C = 0$ for $k \leq n-1$.

We have to show that $f$ has the left lifting property with respect to
any co-$n$-fibration $p$.  By Proposition
\ref{prop:co-n-fibration} and Lemma \ref{lem:n-equivalence}, 
the fiber $F$ of $p$ has the property that $\pi_k F = 0$ for $k \geq n$.  

A lifting problem for $f$ with respect to $p$ has an obstruction
belonging to $[\Omega C, F]$ \cite[Cor.~8.4]{CDI}.  However,
the conditions on the homotopy groups of $\Omega C$ and $F$
guarantee that $[\Omega C, F]$ equals $0$.  
Therefore,
the obstruction to lifting must vanish, and the desired lift
exists.
\end{proof}

We will now show how to build co-$n$-fibrations out of fibrations
whose fibers are Eilenberg-Mac Lane spectra.

\begin{lemma}
\label{lem:co-n-fibration-composition}
Let $m$ and $n$ be any integers.  Any co-$m$-fibration is a retract of a map
that can be factored into a finite
composition of co-$n$-fibrations and fibrations whose fibers
are Eilenberg-Mac Lane spectra.
\end{lemma}

\begin{proof}
If $n \geq m$, then any co-$m$-fibration is a co-$n$-fibration.  Thus,
we may assume that $m > n$.

Let $q: E \map B$ be a co-$m$-fibration.
First use Lemma \ref{lem:n-factor} to factor $q$ into an
$n$-cofibration $j_n: E \map E_n$ followed by a co-$n$-fibration
$q_n:E_n \map B$.  Repeat by factoring $j_{k-1}$ into a $k$-cofibration
$j_k:E \map E_k$ followed by a co-$k$-fibration $q_k:E_k \map E_{k-1}$.
We obtain a diagram
\[
\xymatrix{
E \ar[d]_{j_m} \ar[dr]_{j_{m-1}} \ar[drrr]_{j_n} 
   \ar[drrrr]^{q}     \\
E_m \ar[r]_{q_m} & E_{m-1} \ar[r]_{q_{m-1}} & \cdots \ar[r] & 
E_n \ar[r]_{q_n} & B.             }
\]
Let $p: E_m \map B$ be the composition of the maps $q_n, \ldots, q_m$.
Lifting the $m$-cofibration $j_m$ with respect to the co-$m$-fibration
$q$ shows that $q$ is a retract of $p$.  

The map $q_n$ is a co-$n$-fibration by construction.  
For $n+1 \leq k \leq m$, 
let $F_k$ be the fiber of $q_k$.
Since $q_k$ is a co-$k$-fibration, $\pi_r F_k = 0$
for $r \geq k$ by Lemma \ref{lem:n-equivalence} and Proposition
\ref{prop:co-n-fibration}.
Using that $j_k$ is a $k$-equivalence, that $j_{k-1}$ is a $(k-1)$-equivalence,
and that $j_{k-1} = q_m j_m$, a small diagram chase verifies that
$q_k$ is a $(k-1)$-equivalence.
This implies that $\pi_r F_k = 0$ for $r\leq k-2$.  Hence
$\pi_r F_k$ can only be non-zero when $r = k-1$, so $F_k$ is an 
Eilenberg-Mac Lane spectrum.
\end{proof}

\subsection{Mapping spaces and homotopy classes}
\label{subsctn:mapping}

We next show that $n$-cofibrations interact appropriately with
tensors.  This will
be needed to show that the model structure on pro-spectra is simplicial.
If $X$ is a spectrum and $K$ is a simplicial set, recall that
$X \otimes K$ is defined to be $X \wedge K_+$.

\begin{proposition} \label{prop:n-cofibration-pushout-product}
Suppose that $f: A \map B$ is an $n$-cofibration
and $i:K \map L$ is a cofibration of simplicial sets.
Then the map
\[
g: A \otimes L \amalg_{A \otimes K} B \otimes K
  \map B \otimes L
\]
is also an $n$-cofibration.
\end{proposition}

\begin{proof}
The map $i$ is a transfinite composition of cobase changes of
maps of the form $\bdry{j} \map \Delta[j]$.  
Therefore, the map $g$ is a transfinite composition of cobase changes
of maps of the form 
\[
A \otimes \Delta[j] \amalg_{A \otimes \bdry{j}} B \otimes \bdry{j}
  \map B \otimes \Delta[j].
\]
Since $n$-cofibrations are defined by a left lifting property, 
$n$-cofibrations are preserved by cobase changes and transfinite compositions.
Therefore, we may assume that $i$ is the map $\bdry{j} \map \Delta[j]$.

Spectra  are a simplicial model category
and $f$ is a cofibration by Lemma \ref{lem:n-cofibration-fibration},
so $g$ is also a cofibration.
By
Proposition \ref{prop:n-cofibration}, 
we need only show that $g$ is an $n$-equivalence.
Let $C$ be the cofiber of $f$.  Then the cofiber of $g$ is $C \wedge S^{j}$,
where the simplicial set $S^{j}$ is the sphere $\Delta[j] / \bdry{j}$
based at the image of $\bdry{j}$.
By Lemma \ref{lem:n-equivalence}, we need only show that 
$\pi_{k} (C \wedge S^{j}) = 0$ for all $k \leq n$.
Since $\pi_{k} C = 0$ for all $k \leq n$ by Lemma \ref{lem:n-equivalence}
and $\pi_{k} (C \wedge S^{j}) = \pi_{k-j} C$, it follows that
$\pi_{k} (C \wedge S^{j}) = 0$ for $k \leq j + n$.
This suffices since $j \geq 0$.
\end{proof}

\begin{corollary} \label{cor:n-Map}
Let $A \map B$ be an $n$-cofibration, and let $X \map Y$ be 
a co-$n$-fibration.  The map
\[
f: \Map(B, X) \map \Map(A,X) \times_{\Map(A,Y)} \Map(B,Y)
\]
is an acyclic fibration of simplicial sets.
\end{corollary}

\begin{proof}
This follows from the lifting property characterization of 
acyclic fibrations, adjointness, and 
Proposition \ref{prop:n-cofibration-pushout-product}.
\end{proof}

The next result is a highly technical lemma that will be needed in one place
later.

\begin{lemma}
\label{lem:zero-map}
Let $X \map Y$ be a map of spectra such that 
$* \map X$ is an $(n-1)$-equivalence and such that
the map $\pi_n X \map \pi_n Y$ is zero.
Let
$Z \map *$ be a co-$(n+1)$-equivalence.
Then the map $[Y,Z] \map [X,Z]$ is zero.
\end{lemma}

\begin{proof}
We may assume that $X$, $Y$, and $Z$ are both cofibrant and fibrant.
For each element of $\pi_n X$, choose a representative $S^n \map X$.
This gives a map $\vee S^n \map X$; let $X'$ be its cofiber.  
Note that $\pi_n X'$ equals zero by construction.  Also note that
$\pi_k X'$ is isomorphic to $\pi_k X$ for $k < n$.  
The homotopy groups of $X'$ vanish in dimensions less than or equal to $n$,
and the homotopy groups of $Z$ vanish in dimensions greater than or equal to
$n+1$.  This guarantees
that $[X',Z]$ is zero.  From the exact sequence
\[
[X',Z] \map [X, Z] \map [\vee S^n, Z],
\]
we see that $[X,Z] \map [\vee S^n,Z]$ is injective.

Now the composition $\vee S^n \map X \map Y$ is homotopy trivial
by assumption.  Therefore, the composition
\[
[Y,Z] \map [X, Z] \map [\vee S^n, Z]
\]
is zero.  But the second map is injective, so we can conclude that
the first map is zero.
\end{proof}

\section{Model Structure}\label{sctn:model-structure}

We now define a model structure for pro-spectra.

\begin{definition}\label{defn:weak-equivalence}
A map of pro-spectra $f$ is a \mdfn{$\pi_*$-weak equivalence} if
$f$ is an essentially levelwise $n$-equivalence for every $n$.
\end{definition}

This means that for every $n$, $f$ has a level
representation that is a levelwise an $n$-equivalence.
Beware that $\pi_*$-weak equivalences do not have to be strict
weak equivalences.  The point is that different level representations
may be required for different values of $n$.

The terminology may appear strange at this point.  In Section
\ref{sctn:pro-homotopy}, we will show that the $\pi_*$-weak equivalences
can be recharacterized in terms of of pro-homotopy groups.

\begin{example}
Recall the pro-spectrum $X$ from Section \ref{sctn:K-theory},
where $X_n$ is the wedge $\vee_{k=n}^\infty S^{2k}$.  We claimed in
Section \ref{sctn:K-theory} that the map $* \map X$ is a $\pi_*$-weak
equivalence.  To see that $* \map X$ is an essentially levelwise 
$m$-equivalence, restrict $X$ to the subdiagram $X'$ consisting
of those $X_n$ with $2n > m$.  The subdiagram $X'$ is cofinal in $X$, so 
$X$ and $X'$ are isomorphic as pro-spectra.
The level map from $*$ to $X'$ is a levelwise $m$-equivalence.
\end{example}

The following lemma shows that our model structure is a localization
of the strict model structure (see Section \ref{subsctn:strict}).

\begin{lemma} \label{lem:strict-weak-equivalence}
Strict weak equivalences are $\pi_*$-weak equivalences.
\end{lemma}

\begin{proof}
A levelwise weak equivalence is also a levelwise $n$-equivalence.
\end{proof}

\begin{definition}\label{defn:pro-cofibration}
A map of pro-spectra is a \dfn{cofibration} if it is an essentially
levelwise cofibration.
\end{definition}

\begin{definition}\label{defn:pro-fibration}
A map of pro-spectra is a \dfn{$\pi_*$-fibration} 
if it has the right lifting property
with respect to all $\pi_*$-acyclic cofibrations.
\end{definition}

The terminology here emphasizes that the notion of cofibration is the
same in all known model structures for pro-spectra.
On the other hand,
the fibrations vary among the different model structures.

It requires some work to establish that these definitions
give a model structure on the category of pro-spectra.  We begin by
collecting various technical lemmas.  By the end of this section,
we will be able to prove that the model structure exists.

\subsection{Two-out-of-Three Axiom} \label{sctn:2/3}

This subsection deals with the two-out-of-three
axiom for $\pi_*$-weak equivalences.  Typically this axiom is automatic
from the definition, but we have to do a little work.

\begin{lemma} \label{lem:2/3-n-equivalence}
Suppose that $f$ and $g$ are two composable morphisms of pro-spectra.
If any two of $f$, $g$, and $gf$ are essentially levelwise $n$-equivalences,
then the third is an essentially levelwise $(n-1)$-equivalence.
\end{lemma}

\begin{proof}
The proofs of \cite[Lem.~3.5]{I2} and \cite[Lem.~3.6]{I2}, which concern
the two-out-of-three axiom for essentially levelwise weak equivalences,
can be applied.  To make these proofs work, two formal properties of
$n$-equivalences are required.  First, $n$-equivalences are
preserved by base changes along fibrations and cobase changes along
cofibrations (see Lemma \ref{lem:pullback-pushout-n-equivalence}).
Second, if $f: X \map Y$ and $g: Y \map Z$ are maps of ordinary spectra
and any two of $f$, $g$, and $gf$ are $n$-equivalences, then the third
is an $(n-1)$-equivalence.
\end{proof}

\begin{proposition} \label{prop:2/3}
The $\pi_*$-weak equivalences of pro-spectra given in Definition
\ref{defn:weak-equivalence} satisfy the two-out-of-three 
axiom.
\end{proposition}

\begin{proof}
Let $f$ and $g$ be two composable maps of pro-spectra, and suppose
that two of the maps $f$, $g$, and $gf$ are $\pi_*$-weak equivalences.
By Lemma \ref{lem:2/3-n-equivalence}, 
the third is an essentially levelwise $(n-1)$-equivalence for 
every $n$.  
\end{proof}

\subsection{$\pi_*$-Acyclic Cofibrations} \label{subsctn:acyclic-cofibration}

We shall find it useful to study the essentially levelwise $n$-cofibrations.
Beware that we do not know (yet) that these maps 
are the same as maps that are both essentially 
levelwise cofibrations and essentially levelwise $n$-equivalences.
The difficulty is that the reindexing required to replace a
map by a levelwise cofibration
may not agree with the reindexing required to replace the
same map by a levelwise $n$-equivalence.

\begin{lemma} \label{lem:factor-level-n-cofibration-acyclic-fibration}
Any essentially levelwise $n$-equivalence factors into 
an essentially levelwise $n$-cofibration followed by a
strict acyclic fibration.
\end{lemma}

\begin{proof}
We may assume that $f$ is a level map that is a levelwise
$n$-equi\-val\-ence.  Use the method of \cite[Lem.~4.6]{I2} to factor $f$ into
a levelwise cofibration $i$ followed by a strict acyclic
fibration $p$.  By \cite[Lem.~4.4]{I2}, $p$ is also a levelwise
acyclic fibration.

For each $s$, we have $f_s = p_s i_s$.  Since $f_s$ is an $n$-equivalence
and $p_s$ is a weak equivalence, it follows that $i_s$ is also an
$n$-equivalence.
Now use Proposition \ref{prop:n-cofibration} to conclude that 
$i$ is a levelwise $n$-cofibration.
\end{proof}

\begin{proposition} \label{prop:level-n-cofibration}
A map is a cofibration and essentially levelwise $n$-equivalence if
and only if it is an essentially levelwise $n$-cofibration.
\end{proposition}

\begin{proof}
First suppose that $i$ is an essentially levelwise $n$-cofibration.
Then $i$ is an essentially levelwise cofibration because every
$n$-cofibration is a cofibration.  Similarly, $i$ is an essentially
levelwise $n$-equivalence because every $n$-cofibration is an
$n$-equivalence.

For the other direction,
let $i$ be a cofibration and essentially levelwise $n$-equi\-valence.
By Lemma \ref{lem:factor-level-n-cofibration-acyclic-fibration}, 
$i$ factors into an essentially
levelwise $n$-cofibration $j$ followed by a strict acyclic fibration $p$.
Then $p$ has the right lifting property with respect to the cofibration
$i$ because of the strict structure, 
so $i$ is a retract of $j$.  Essentially levelwise $n$-cofibrations
are closed under retract by \cite[Cor.~5.6]{I3}.
\end{proof}

\subsection{$\pi_*$-Fibrations} \label{subsctn:fibrations}

\begin{lemma} \label{lem:strict-fibration}
Every $\pi_*$-fibration is a strict fibration,
and every strict acyclic fibration is a $\pi_*$-acyclic fibration.
\end{lemma}

\begin{proof}
For the first claim, observe that
Lemma \ref{lem:strict-weak-equivalence} guarantees  every strict
acyclic cofibration is a $\pi_*$-acyclic cofibration.  Now use the
lifting property definitions of $\pi_*$-fibrations and strict fibrations.

For the second claim, recall that strict acyclic fibrations have the
right lifting property with respect to all cofibrations and therefore
with respect to $\pi_*$-acyclic cofibrations.  This means that a
strict acyclic fibration is a $\pi_*$-fibration.  To show that
it is also a $\pi_*$-weak equivalence, use 
Lemma \ref{lem:strict-weak-equivalence}.
\end{proof}

\begin{lemma} \label{lem:levelwise-fibration}
Every $\pi_*$-fibration is an essentially levelwise fibration.
\end{lemma}

\begin{proof}
This follows from Lemma \ref{lem:strict-fibration} and 
the fact that strict fibrations are essentially levelwise fibrations
\cite[Lem.~4.5]{I2}.
\end{proof}

Next we produce some examples of $\pi_*$-fibrations.

\begin{lemma} \label{lem:special-fibration}
Let $X \map Y$ be a co-$m$-fibration for some $m$.  Then the constant
map $p:cX \map cY$ is a $\pi_*$-fibration.
\end{lemma}

\begin{proof}
We show that $p$ has the desired right lifting property.  Let $i: A \map B$
be a $\pi_*$-acyclic cofibration, so $i$ is an essentially levelwise
$m$-equivalence.
By Proposition \ref{prop:level-n-cofibration}, 
we may assume
that $i$ is a levelwise $m$-cofibration.

Suppose given a square
\[
\xymatrix{
A \ar[r] \ar[d]_i & cX \ar[d]^p \\
B \ar[r] & cY                  }
\]
of pro-spectra.  This square is represented by a square
\[
\xymatrix{
A_s \ar[r] \ar[d]_{i_s} & X \ar[d] \\
B_s \ar[r] & Y                  }
\]
of spectra for some $s$.  Now $i_s$ is an $m$-cofibration
and $X \map Y$ is a co-$m$-fibration, so this last square
has a lift.  The lift represents the desired lift.
\end{proof}

\subsection{Small object argument} 
\label{subsctn:fibration}

Eventually we will produce factorziations with a dual version of the
generalized small object argument \cite{C}.
The next results are the technical details that allow us to apply
this technique.

\begin{definition}
\label{defn:F}
Given a level map $f:X \map Y$, let \mdfn{$F(f)$} be the
set of fibrations of spectra defined as follows.  For each $s$ and
each $n$, consider the functorial factorization of $f_s: X_s \map Y_s$
into an $n$-cofibration $i_{s,n}:X_s \map Z_{s,n}$ followed by a
co-$n$-fibration $p_{s,n}: Z_{s,n} \map Y_s$ as in
Lemma \ref{lem:n-factor}.  Let $F(f)$ be the set of 
all such maps $p_{s,n}$.
\end{definition}

\begin{lemma}
\label{lem:detect-acyclic-cofibrations}
A map $i: A \map B$ is a $\pi_*$-acyclic cofibration if and only if
it has the left lifting property with respect to all constant pro-maps
$cX \map cY$ in which $X \map Y$ is a co-$m$-fibration for some $m$.
\end{lemma}

\begin{proof}
One implication is shown in Lemma \ref{lem:special-fibration}.
For the other implication, suppose that $i$ has the desired lifting
property.  Since acyclic fibrations of spectra are co-$m$-fibrations,
$i$ has the left lifting property with respect to all maps
$cX \map cY$ in which $X \map Y$ is an acyclic fibration of spectra.
By \cite[Prop.~5.5]{I2}, this implies that $i$ is a cofibration.

Fix an $n$.
We show that $i$ is an essentially levelwise $n$-equivalence.
From the previous paragraph, 
we may assume that $i$ is a levelwise cofibration.  

Consider the square
\[
\xymatrix{
X \ar[r]\ar[d] & cZ_{s,n} \ar[d] \\
Y \ar[r] & cY_s            }
\]
of pro-spectra, where the map $X \map cZ_{s,n}$ is the 
composition of the canoncial map $X \map cX_s$ together with the
map $i_{s,n}:X_s \map Z_{s,n}$ (see Definition \ref{defn:F})
and the
map $Y \map cY_s$ is the canonical
map.
Our assumption gives us a lift in this diagram because $Z_{s,n} \map Y_s$
is a co-$n$-fibration.
This means that we have a diagram
\[
\xymatrix{
X_t \ar[r] \ar[d] & X_s \ar[r] & Z_{s,n} \ar[d] \\
Y_t \ar[urr]\ar[r] & Y_s \ar[r] & Y_s   }
\]
for some $t$, which can be rewritten as
\[
\xymatrix{
X_t \ar[r]\ar[d] & X_s \ar[r]\ar[d] & X_s \ar[d] \\
Y_t \ar[r] & Z_{s,n} \ar[r] & Y_s.                }
\]

Finally, Lemma \ref{lem:pro-isomorphism} shows that
the objects $Z_{s,n}$ can be assembled into a pro-spectrum that
is isomorphic to $Y$.
Thus, the maps
$X_s \map Z_{s,n}$ give a level representation of $f$.  
Each map 
$X_s \map Z_{s,n}$ is a levelwise $n$-equivalence, so
$X \map Y$ is an essentially levelwise $n$-equivalence.
\end{proof}

\begin{lemma}
\label{lem:F}
Consider a square
\[
\xymatrix{
X \ar[r]\ar[d]_f & cE \ar[d]^{p} \\
Y \ar[r] & cB                    }
\]
of pro-spectra in which 
$p$ is a constant pro-map such that $E \map B$ is a 
co-$n$-fibration for some $n$.  
This diagram factors as 
\[
\xymatrix{
X \ar[r]\ar[d]_f & cZ_{s,n} \ar[r]\ar[d] & cE \ar[d]^{p} \\
Y \ar[r] & cY_s \ar[r] & cB                    }
\]
for some $p_{s,n}:Z_{s,n} \map Y$ belonging to $F(f)$.
\end{lemma}

\begin{proof}
We may assume that $f$ is a level map.
The original square is represented by a diagram
\[
\xymatrix{
X_s \ar[r]\ar[d] & E \ar[d] \\
Y_s \ar[r] & B                    }
\]
for some $s$.
This gives us a square
\[
\xymatrix{
X_s \ar[r]\ar[d]_{i_{s,n}} & E \ar[d]  \\
Z_{s,n} \ar[r] & B                     }
\]
in which the bottom horizontal map is the composition of $p_{s,n}$ 
with the given map $Y_s \map B$.
Note that the left vertical map is an $n$-cofibration
(see Definition \ref{defn:F}) and the right vertical map
is a co-$n$-fibration.  Therefore, a lift $h$ exists in this diagram.

Such a lift $h$ gives us a diagram
\[
\xymatrix{
X_s \ar[r]\ar[d]_{f_s} & Z_{s,n} \ar[r]^h \ar[d]_{p_{s,n}} & E \ar[d] \\
Y_s \ar[r]_= & Y_s \ar[r] & B,                                        }
\]
and this produces the desired factorization.
\end{proof}

\subsection{The $\pi_*$-model structure}
\label{subsctn:model-structure}

We are now ready to prove that the model structure axioms are satisfied.

\begin{theorem}
\label{thm:pro-model-structure}
The cofibrations, $\pi_*$-weak equivalences, and
$\pi_*$-fibrations are a simplicial proper
model structure on the category of pro-spectra.
\end{theorem}

We call this the \mdfn{$\pi_*$-model structure} for pro-spectra.

\begin{proof}
The category of pro-spectra
has all limits and colimits since the category of spectra does
\cite[Prop.~11.1]{I1}.
The two-out-of-three axiom for 
$\pi_*$-weak equivalences is not automatic;  we
proved this in Proposition \ref{prop:2/3}. 
Retracts preserve essentially levelwise properties \cite[Cor.~5.6]{I3}.
Therefore, retracts preserve cofibrations and $\pi_*$-weak equivalences.
Retracts preserve $\pi_*$-fibrations because retracts 
preserve lifting properties.

See \cite[Lem.~4.6]{I2} for factorizations into 
cofibrations followed by maps that are strict acyclic fibrations.
By Lemma \ref{lem:strict-fibration}, strict
acyclic fibrations are $\pi_*$-acyclic fibrations.
This gives factorizations into cofibrations followed by $\pi_*$-acyclic 
fibrations.

We next construct factorizations into $\pi_*$-acyclic cofibrations followed by
$\pi_*$-fi\-bra\-tions.
The generalized small object argument \cite{C} can be applied to the class
of maps $cX \map cY$ such that $X \map Y$ is a co-$m$-fibration
for some $m$.  Actually, we are applying the categorical dual.
The cosmallness hypothesis is proved in \cite[Prop.~3.3]{CI}.
The other hypothesis is Lemma \ref{lem:F}.

We use Lemma \ref{lem:detect-acyclic-cofibrations} to conclude that
the first map in the factorization is a $\pi_*$-acyclic cofibration.
To conclude that the second map is a $\pi_*$-fibration, we 
use Lemma \ref{lem:special-fibration} and note that
the second map 
is constructed as a composition of a transfinite tower of maps
that are base changes of 
maps of the form
$cX \map cY$ such that $X \map Y$ is a co-$m$-fibration
for some $m$.  Now apply the formal properties of right lifting properties.

One of the lifting axioms follows by definition.
The other follows from the retract argument \cite[Prop.~7.2.2]{Hi}.
In more detail,
any $\pi_*$-acyclic fibration $p$ can be factored
into a cofibration $i$ followed by a strict acyclic fibration $p$.
Then $i$ is a $\pi_*$-weak equivalence by the two-out-of-three
axiom and the fact that strict weak equivalences are
$\pi_*$-weak equivalences (see Lemma \ref{lem:strict-weak-equivalence}).  
Hence $p$ has the right lifting property with respect to $i$,
so $p$ is a retract of $q$.  It follows that $p$ is a strict
acyclic fibration,
so it has the right lifting property with respect to all cofibrations.

The simplicial structure is analogous to the simplicial 
structure for pro-spaces \cite[\S~16]{I1}.
Beware that the definitions of tensor and cotensor are 
straightforward for finite simplicial sets but are slightly subtle
in general.
We need to show that if $i:K \map L$ is a cofibration of 
finite simplicial
sets and $f: A \map B$ is a cofibration of pro-spectra, then the map
\[
g: A \otimes L \coprod_{A \otimes K} B \otimes K \map B \otimes L
\]
is a cofibration of pro-spectra that is a $\pi_*$-weak equivalence
if either $i$ is an acyclic cofibration or 
$f$ is a $\pi_*$-acyclic cofibration.
The fact that $g$ is a cofibration follows from the fact that the
strict model structure is simplicial \cite[Thm.~4.16]{I2}.  The case
when $i$ is an acyclic cofibration also follows from the strict
structure.  

It remains to assume that $f$ is a $\pi_*$-acyclic cofibration.  Given
any $n$, we may assume that $f$ is a levelwise $n$-cofibration
by Proposition \ref{prop:level-n-cofibration}.  Since tensors with
finite simplicial sets can be constructed levelwise and since pushouts
can be constructed levelwise \cite[App.~4.2]{AM}, it follows from 
Proposition \ref{prop:n-cofibration-pushout-product} that $g$ is 
a levelwise $n$-cofibration.  This means that $g$ is a levelwise
$n$-equivalence for every $n$, so $g$ is a $\pi_*$-weak equivalence.

For right properness,
consider a pullback square
\[
\xymatrix{
W \ar[r]^q \ar[d]_g & X \ar[d]^f \\
Y \ar[r]_p & Z               }
\]
in which $f$ is a $\pi_*$-weak equivalence and $p$ is a $\pi_*$-fibration.
We want to show that $g$ is also a $\pi_*$-weak equivalence.
Lemma \ref{lem:strict-fibration} implies that $p$ is a strict fibration.
Therefore, the proof of \cite[Thm.~4.13]{I2} can be applied to show
that base changes of essentially levelwise $n$-equivalences
along $\pi_*$-fibrations
are again essentially levelwise $n$-equivalences.
We need Lemma \ref{lem:pullback-pushout-n-equivalence} for the proof
to work.

The proof of left properness is dual.
\end{proof}

\begin{remark}
The proof of properness for pro-spaces given in \cite[Prop.~17.1]{I1}
is incorrect, but the techniques of \cite[Thm.~4.13]{I2} can be
used to fix it.
\end{remark}

We write $\Map(X, Y)$ 
for the simplicial mapping space of pro-maps
from $X$ to $Y$.  More precisely, we have the formula
\[
\Map(X, Y)  = \lim_{s} \colim_{t} \Map(X_{t}, Y_{s}).
\]

Constructing cofibrant replacements is straightforward.  Given a 
pro-spectrum $X$, we just take a levelwise cofibrant replacement.
In Section \ref{sctn:fibrant}, we will show that
constructing $\pi_*$-fibrant replacements is a bit more complicated.
Let $X$ be a pro-spectrum indexed by a cofiltered category $I$. 
Define a new pro-spectrum $PX$ indexed by $I \times \Z$ as follows.
For every pair $(s,n)$, let $PX_{(s,n)} = P_{n} X_{s}$ be the
$n$th Postnikov section of $X_{s}$.  Finally, take a strict fibrant
replacement for $PX$.  The resulting pro-spectrum is a $\pi_*$-fibrant
replacement for $X$.

\subsection{Stable model structure}
\label{subsctn:stable}

Recall that the functors $- \wedge S^1$ and $\Map(S^1,-)$ are
defined levelwise for pro-spectra.  
In this section, we will show that the $\pi_*$-model structure is 
stable in the sense that these functors are a Quillen equivalence
from the $\pi_*$-model structure to itself.

\begin{lemma}
\label{lem:pi-Quillen-pair}
The functors $- \wedge S^1$ and $\Map(S^1,-)$ are a Quillen adjoint pair
from the $\pi_*$-model structure on pro-spectra to itself.
\end{lemma}

\begin{proof}
On spectra, $- \wedge S^1$ preserves cofibrations.  Therefore, it
preserves levelwise cofibrations and thus essentially levelwise 
cofibrations.  

On spectra, $- \wedge S^1$ takes $n$-cofibrations to $(n+1)$-cofibrations.
With the help of 
Proposition \ref{prop:level-n-cofibration}, this shows that $- \wedge S^1$
preserves $\pi_*$-acyclic cofibrations. 
\end{proof}

\begin{lemma}
\label{lem:pi-suspension}
Let $A$ and $B$ be cofibrant pro-spectra.  A map $f: A \map B$ is
a $\pi_*$-weak equivalence if and only if $f \wedge S^1$ is a 
$\pi_*$-weak equivalence.
\end{lemma}

\begin{proof}
To simplify notation, write $F$ for the functor $- \wedge S^1$
and $G$ for the functor $\Map(S^1,-)$.

One direction follows immediately from Lemma \ref{lem:pi-Quillen-pair}
and the fact that left Quillen functors preserve weak equivalences between
cofibrant objects \cite[Prop.~8.5.7]{Hi}.

For the other direction, suppose that $Ff$ is a $\pi_*$-weak equivalence.
Factor $f$ into a cofibration $i: A \map C$ followed by a strict
acyclic fibration $p: C \map B$.  
The map $Fp$ is a strict weak equivalence
because left Quillen functors preserve weak equivalences between
cofibrant objects.  By the two-out-of-three axiom, $Fi$ is also
a $\pi_*$-weak equivalence.  
By the two-out-of-three axiom again,
$f$ is a $\pi_*$-weak equivalence if $i$ is a 
$\pi_*$-weak equivalence.
Therefore, it suffices to show that the cofibration $i$ is a 
$\pi_*$-weak equivalence.

We use the 
lifting characterization of Lemma \ref{lem:detect-acyclic-cofibrations}
to show that $i$ is a $\pi_*$-acyclic cofibration.
We may assume that
$A$ and $C$ are both levelwise cofibrant.  Using the level replacement of
\cite[App.~3.2]{AM}, 
we may reindex $A$ and $C$ in such a
way that $i$ is a level map and 
$A$ and $C$ are still levelwise cofibrant.
However, we are not allowed to assume that $i$ is
a levelwise cofibration because this may require a different reindexing.

Suppose given a lifting problem
\[
\xymatrix{
A \ar[r]\ar[d] & cX \ar[d] \\
C \ar[r] & cY,             }
\]
where $X \map Y$ is a co-$m$-fibration for some $m$.  This diagram
of pro-spectra is represented by a diagram
\[
\xymatrix{
A_s \ar[r]\ar[d] & X \ar[d] \\
C_s \ar[r] & Y             }
\]
of spectra for some $s$.  First factor, $C_s \map Y$ into a
cofibration $C_s \map \tilde{Y}$ followed by an acyclic fibration
$\tilde{Y} \map Y$.  Then factor the map $A_s \map \tilde{Y} \times_Y X$
into a cofibration $A_s \map \tilde{X}$ followed by an acyclic fibration
$\tilde{X} \map \tilde{Y} \times_Y X$.  This gives us a diagram
\[
\xymatrix{
A_s \ar[r]\ar[d] & \tilde{X} \ar[r]\ar[d] & X \ar[d] \\
C_s \ar[r] & \tilde{Y} \ar[r] & Y             }
\]
in which the map $\tilde{X} \map \tilde{Y}$ is a fibration, the maps
$\tilde{X} \map X$ and $\tilde{Y} \map Y$ are weak equivalences,
and the spectra $\tilde{X}$ and $\tilde{Y}$ are cofibrant.
By Proposition \ref{prop:co-n-fibration}, $\tilde{X} \map \tilde{Y}$ is also a
co-$m$-fibration.
Now we have a diagram
\[
\xymatrix{
A \ar[r]\ar[d] & c\tilde{X} \ar[r]\ar[d] & cX \ar[d] \\
C \ar[r] & c\tilde{Y} \ar[r] & cY             }
\]
of pro-spectra.  We want to show that the outer rectangle has a lift,
so it suffices to show that the left square has a lift.

Let $\hat{F} \tilde{Y}$ be a fibrant replacement for $F\tilde{Y}$.
Factor the composition $F\tilde{X} \map \hat{F}\tilde{Y}$ into
an acyclic cofibration $F\tilde{X} \map \hat{F} \tilde{X}$ followed
by a fibration $\hat{F} \tilde{X} \map \hat{F} \tilde{Y}$.  Note
that $\hat{F} \tilde{X}$ is a fibrant replacement for $F\tilde{X}$.

The maps $\tilde{X} \map G\hat{F} \tilde{X}$ and
$\tilde{Y} \map G\hat{F} \tilde{Y}$ are weak equivalences
because $F$ and $G$ are a Quillen equivalence on spectra.  Here is
where we use that $\tilde{X}$ and $\tilde{Y}$ are cofibrant.
Because $G$ is a right Quillen functor, the map 
$G\hat{F} \tilde{X} \map G\hat{F} \tilde{Y}$ is also a fibration.
Moreover, this fibration is a co-$m$-fibration by Proposition
\ref{prop:co-n-fibration}.

Now consider the diagram
\[
\xymatrix{
A \ar[r]\ar[d] & c\tilde{X} \ar[r]\ar[d] & cG\hat{F} \tilde{X} \ar[d] \\
C \ar[r] & c\tilde{Y} \ar[r] & cG\hat{F} \tilde{Y}                    }
\]
of pro-spectra.  A lift exists in the outer rectangle by adjointness
and Lemma \ref{lem:detect-acyclic-cofibrations} applied to
the $\pi_*$-acyclic cofibration $Fi: FA \map FC$.
By \cite[Prop.~3.2]{CDI}, a lift exists in the left square also. 
\end{proof}

\begin{theorem}
\label{thm:pro-stable}
The functors $- \wedge S^1$ and $\Map(S^1, -)$ are a Quillen equivalence 
from the $\pi_*$-model structure on pro-spectra to itself.
\end{theorem}

\begin{proof}
As before, to simplify the notation, write $F$ for $- \wedge S^1$
and $G$ for $\Map(S^1,-)$.

Suppose that $g: X \map GY$ is any map such that $X$ is cofibrant
and $Y$ is $\pi_*$-fibrant.
We want to show that $g$ is a 
$\pi_*$-weak equivalence if and only if its adjoint
$f: FX \map Y$ is a $\pi_*$-weak equivalence.

Factor $g$ into a cofibration $i:X \map Z$ followed by a strict
acyclic fibration $p: Z \map GY$.
The adjoint $p':FZ \map Y$ is a strict weak equivalence
because $F$ and $G$ are a Quillen equivalence on the strict
model structure as shown in Theorem \ref{thm:strict-stable}.  Here we are using
that $Y$ is strict fibrant by Lemma \ref{lem:strict-fibration}.

The adjoint $f$ is the composition of $Fi$ with $p'$.  By the
two-out-of-three axiom, $f$ is a $\pi_*$-weak equivalence if and only
if $Li$ is a $\pi_*$-weak equivalence.  Because $X$ and $Z$ are
cofibrant, Lemma \ref{lem:pi-suspension} tells us that $Fi$ is a 
$\pi_*$-weak equivalence if and only if $i$ is a $\pi_*$-weak equivalence.
Finally, the two-out-of-three axiom implies that $i$ is a 
$\pi_*$-weak equivalence if and only if $g$ is a $\pi_*$-weak equivalence.
\end{proof}

\section{$\pi_*$-fibrant pro-spectra}
\label{sctn:fibrant}

\begin{theorem}
\label{thm:fibrant}
A pro-spectrum $X$ is $\pi_*$-fibrant if and only if it
is strict fibrant and essentially levelwise
fibrant and bounded above.
\end{theorem}

That $X$ is levelwise fibrant and bounded above 
means that each $X_s$ is fibrant and bounded above 
(see Definition \ref{defn:bounded});
we require no uniformity on the dimension in which the homotopy groups
vanish.

The following proof is similar to the proof of \cite[Prop.~4.9]{CI}.

\begin{proof}
First suppose that $X$ is strict fibrant and essentially levelwise
fibrant and bounded above.
We will show that for every $\pi_*$-acyclic cofibration $i: A \map B$,
the map $f:\Map(B,X) \map \Map(A,X)$ of simplicial sets is an acyclic
fibration.  By the usual adjointness arguments, this will show that
$X \map *$ has the desired lifting property.  Since
$X$ is strict fibrant and $i$ is a cofibration, we already know that
$f$ is a fibration.  It remains to show that $f$ is a weak equivalence.

We may assume that each $X_s$ is fibrant and bounded above.
We showed in Lemma \ref{lem:special-fibration} that the constant
pro-spectrum $cX_s$ is $\pi_*$-fibrant.  
Therefore the map $\Map(B, cX_s) \map \Map(A, cX_s)$ is an acyclic fibration
and in particular a weak equivalence.

Since $X$ is levelwise fibrant and strict fibrant,
Proposition \ref{prop:strict-map} implies that 
the mapping space $\Map(A, X)$ 
is weakly equivalent to $\holim_s \Map(A, cX_s)$
(and similarly for $\Map(B, X)$).
Homotopy limits preserve weak equivalences,
so we conclude that $\Map(B,X) \map \Map(A,X)$ is a weak equivalence.
This completes one implication.

Now suppose that $X$ is $\pi_*$-fibrant.  Then $X$ is strict fibrant
by Lemma \ref{lem:strict-fibration}.
It remains to show that $X$ is essentially
levelwise fibrant and bounded above.

Consider the factorization $X \map Y \map *$ of the map
$X \map *$ into a $\pi_*$-acyclic cofibration followed by a 
$\pi_*$-fibration by means of the generalized small object argument
(see Section \ref{sctn:model-structure}).
Now $X$ is a retract of $Y$ because $X$ is $\pi_*$-fibrant
and $X \map Y$ is a $\pi_*$-acyclic cofibration.
The class of pro-objects having any 
property essentially levelwise is closed under retracts 
\cite[Thm.~5.5]{I3},
so it suffices to consider $Y$.

Recall that $Y \map *$ is constructed as a composition of a transfinite tower
\[
\cdots \map Y_\beta \map \cdots \map Y_2 \map Y_1 \map *,
\]
where each map $Y_{\beta+1} \map Y_\beta$ is a base change of a product
of maps of the form $cE \map cB$ with $E \map B$ a co-$m$-fibration
for some $m$.
The class of pro-objects having any 
property essentially levelwise is closed under cofiltered limits
\cite[Thm.~5.1]{I3},
so it suffices to consider each $Y_{\beta}$.

We proceed by transfinite induction.  When $\beta$ is a limit ordinal,
\cite[Thm.~5.1]{I3} again tells us that $Y_\beta$ 
is essentially levelwise fibrant and bounded above.

It only remains to consider the case when $\beta$ is a successor ordinal.
We are assuming that $Y_{\beta-1}$
is levelwise fibrant and bounded above.
We may take a level representation for the diagram
\[
\xymatrix{
Y_{\beta-1} \ar[r] & \prod_a cB_a & \prod_a cE_a, \ar[l] }
\]
where each $E_a \map B_a$ is a co-$m$-fibration for some $m$.
Note that $m$ depends on $a$.
We construct $Y_{\beta}$ by taking the levelwise fiber product.
It is possible to construct a level representation for the above 
diagram in such 
a way that the replacement for $Y_{\beta-1}$ is a diagram of objects
that already appeared in the original $Y_{\beta-1}$.  This means
that the new $Y_{\beta-1}$ is still levelwise fibrant and bounded above.

The construction of arbitrary products in pro-categories
\cite[Prop.~11.1]{I1} shows that the map
$\prod_a cE_a \map \prod_a cB_a$
is levelwise a finite product of maps of the form
$E_a \map B_a$.
A finite product of maps that are co-$m$-fibrations for some $m$
is again a co-$m$-fibration for some $m$, so the map 
$Y_\beta \map Y_{\beta-1}$ is levelwise 
a base change of a co-$m$-fibration for some $m$.
Since co-$m$-fibrations are closed under base change,
we conclude that $Y_\beta \map Y_{\beta-1}$ is a levelwise
co-$m$-fibration.
It follows immediately that $Y_{\beta}$ is levelwise fibrant
and bounded above.
\end{proof}

\begin{remark}
\label{rem:ML-fibration}
Similarly to \cite[Prop.~6.6]{I1} and \cite[Defn.~4.2]{I2},
it is possible to give a concrete description of the
$\pi_*$-fibrations. 
Recall that a directed set is cofinite if for every $s$,
there are only finitely many $t$ such that $t \leq s$.  Suppose that
$f:X \map Y$ is a level map indexed by a cofinite directed set such that
each map
\[
X_s \map Y_s \times_{\lim_{t<s} Y_t} \lim_{t<s} X_t
\]
is a co-$m$-fibration for some $m$.  Here $m$ depends on $s$.  
Then $f$ is a $\pi_*$-fibration.  Up to retract, every $\pi_*$-fibration
is of this form.
\end{remark}

The following corollary simplifies the construction of $\pi_*$-fibrant
replacements.

\begin{corollary}
\label{cor:fibrant}
If $Y$ is an essentially levelwise bounded above pro-spectrum,
then there is a strict fibrant replacement $\hat{Y}$ for $Y$
such that $\hat{Y}$ is also a $\pi_*$-fibrant replacement for $Y$.
\end{corollary}

\begin{proof}
We may assume that $Y$ is levelwise bounded above.
Factor the map $Y \map *$ into a strict acyclic cofibration
$Y \map \hat{Y}$ followed by a strict fibration $\hat{Y} \map *$
using the method of \cite[Lem.~4.7]{I2}.  This particular construction
gives that $Y \map \hat{Y}$ is a levelwise weak equivalence
and $\hat{Y} \map *$ is a levelwise fibration;
thus $\hat{Y}$ is levelwise fibrant and bounded above.
Now Theorem \ref{thm:fibrant} implies that $\hat{Y}$ is
$\pi_*$-fibrant.
\end{proof}

The next corollary simplifies the computation of mapping spaces
of pro-spectra.

\begin{corollary}
\label{cor:fibrant2}
Let $X$ be a cofibrant pro-spectrum, and let $Y$ be a levelwise
fibrant bounded above pro-spectrum with $\pi_*$-fibrant replacement $\hat{Y}$.
Then the homotopically correct mapping space
$\Map(X, \hat{Y})$ is weakly equivalent to
$\holim_s \colim_t \Map(X_t, Y_s)$.
\end{corollary}

\begin{proof}
Because the mapping space is homotopically correct, it doesn't matter
which $\pi_*$-fibrant replacement $\hat{Y}$ we consider.  Thus, we
may take the one from Corollary \ref{cor:fibrant}.
Because $\hat{Y}$ is a strict fibrant replacement for $Y$,
Proposition \ref{prop:strict-map} can be applied.
\end{proof}

\section{Homotopy classes of maps of pro-spectra}
\label{sctn:homotopy-class}

Let \mdfn{$[X, Y]_{\pros}$} be the set of weak homotopy classes from $X$ to $Y$
in the $\pi_*$-homotopy category of pro-spectra.
Let
\mdfn{$[X, Y]_{\pros}^{r}$} be the set of 
weak homotopy classes of degree $r$ from $X$ to $Y$.  
For all $r$, $[X, Y]_{\pros}^{r}$ is equal to 
\[
[\Sigma^{-r} X, Y]_{\pros} = [X, \Sigma^{r} Y]_{\pros},
\]
where $\Sigma^r$ equals $\Omega^{-r}$ if $r < 0$.

The mapping space $\Map(X,Y)$
is related to homotopy classes in the following way.
For every cofibrant $X$, fibrant $Y$, and $r \geq 0$,
\[
[X, Y]^{-r}_{\pros} \cong \pi_{r} \Map(X, Y) \cong
\pi_{0} \Map(\Sigma^{r} X, Y).
\]

\begin{proposition} \label{prop:into-constant}
Let $X$ be a pro-spectrum and $Y$ be a bounded above spectrum.
Then $[X, cY]^{r}_{\pros}$ 
is equal to $\colim_{s} [X_{s},Y]^{r}$.
\end{proposition}

\begin{proof}
We may assume that $X$ is levelwise cofibrant and that
$\Sigma^r Y$ is a fibrant spectrum.
We must calculate homotopy classes of maps from $X$ to $c\Sigma^{r} Y$.
Now $\Sigma^r Y$ is 
bounded above since its homotopy groups are just the shifted
homotopy groups of $Y$.
Thus Theorem \ref{thm:fibrant} 
tells us that
the constant pro-spectrum $c\Sigma^{r}Y$ is already $\pi_*$-fibrant.
Therefore,
\[
[X, cY]^{r}_{\pros} \cong
\pi_{0} \Map(X, c\Sigma^{r} Y) \cong
\colim_{s} \pi_{0} \Map(X_{s}, \Sigma^{r} Y) \cong
\colim_{s} [X_{s}, Y ]^{r}.
\]
\end{proof}

Proposition \ref{prop:into-constant} is certainly false if $Y$ is not
bounded above.  For example, let $X$ be the pro-spectrum from
Section \ref{sctn:K-theory}, and let $Y$ be the spectrum $KU$.
Since $X$ is contractible, $[X, cKU]_{\pros}$ is zero.  On the other
hand, we showed in Section \ref{sctn:K-theory} that
$\colim_s [X_s, KU]$ is uncountable.

\begin{lemma}
\label{lem:pro-connected-Postnikov}
Let $* \map X$ be an essentially levelwise $n$-cofibration,
and let $Y \map *$ be an essentially levelwise co-$n$-fibration.
Then the homotopically correct mapping space 
$\Map(X, \hat{Y})_{\pros}$ is trivial, where $\hat{Y}$ is a
$\pi_*$-fibrant replacement.
\end{lemma}

\begin{proof}
We may assume that $* \map X$ is a levelwise $n$-cofibration
and that $Y \map *$ is a levelwise co-$n$-fibration.
In particular, this implies that $Y$ is levelwise bounded above.

By Corollary \ref{cor:n-Map}, each space
$\Map(X_t, Y_s)$ is contractible.  
Therefore, the filtered colimit
$\colim_t \Map(X_t, Y_s)$ is also contractible.
It follows that the cofiltered homotopy limit 
$\holim_s \colim_t \Map(X_t, Y_s)$ is still contractible.
Finally, Corollary \ref{cor:fibrant2} implies that this homotopy
limit is weakly equivalent to $\Map(X, \hat{Y})$.
\end{proof}

\begin{corollary}
\label{cor:pro-connected-Postnikov}
Let $* \map X$ be an essentially levelwise $n$-equivalence,
and let $Y \map *$ be an essentially levelwise co-$n$-equivalence.
Then $[X, Y]_{\pros}$ is zero.
\end{corollary}

\begin{proof}
We may assume that $X$ is cofibrant, so Proposition 
\ref{prop:level-n-cofibration}
implies that we may assume that $* \map X$ is a levelwise $n$-cofibration.

We may assume that that $Y \map *$ is a levelwise co-$n$-equivalence.
By taking a levelwise fibrant replacement, we may further assume that
$Y \map *$ is a levelwise fibration; thus $Y \map *$ is a levelwise
co-$n$-fibration by Proposition \ref{prop:co-n-fibration}.

Now the hypotheses of Lemma \ref{lem:pro-connected-Postnikov} are
satisfied, so the homotopically correct mapping space is contractible.
This implies that $[X,Y]_{\pros}$ is trivial.
\end{proof}

\section{Pro-homotopy groups}
\label{sctn:pro-homotopy}

In this section, we give an alternative characterization of the $\pi_*$-weak
equivalences of Definition \ref{defn:weak-equivalence}.
First we must discuss the stable homotopy pro-groups of a pro-spectrum.
Since $\pi_{k}$ is a functor on spectra, we may apply it 
objectwise to any pro-spectrum $X$ to obtain a pro-group $\pi_{k} X$.

\begin{proposition}\label{prop:cofiber-sequence}
Let $i: A \map B$ be a cofibration with cofiber $C$.
Then there is a long exact sequence
\[
\cdots \map \pi_{k} A \map \pi_{k} B \map \pi_{k} C \map \pi_{k-1} A \cdots
\]
of pro-homotopy groups.
\end{proposition}

To understand what exactness means for this sequence,
see \cite[App.~4.5]{AM}
for a discussion of the abelian structure on the
category of pro-abelian groups.

\begin{proof}
We may suppose that $i$ is a level cofibration, and we may construct
$C$ as the levelwise cofiber of $i$ because finite colimits
in pro-categories can be constructed levelwise \cite[App.~4.2]{AM}.
Now for every $s$, we have a long exact sequence
\[
\cdots \map \pi_{k} A_{s} \map \pi_{k} B_{s} \map \pi_{k} C_{s} \map
\pi_{k-1} A_{s} \map \cdots
\]
of abelian groups.  These sequences assemble to give the 
desired sequence.
\end{proof}

\begin{proposition}\label{prop:fiber-sequence}
Let $p: X \map Y$ be a $\pi_*$-fibration with fiber $F$.
Then there is a long exact sequence
\[
\cdots \map \pi_{k} F \map \pi_{k} X \map \pi_{k} Y \map \pi_{k-1} F \cdots
\]
of pro-homotopy groups.
\end{proposition}

\begin{proof}
By Lemma \ref{lem:levelwise-fibration},
we may assume that
$p$ is a levelwise fibration.
We may construct
$F$ as the levelwise fiber of $p$ because finite limits
in pro-categories can be constructed levelwise \cite[App.~4.2]{AM}.
Now for every $s$, we have a long exact sequence
\[
\cdots \map \pi_{k} F_{s} \map \pi_{k} X_{s} \map \pi_{k} Y_{s} \map
\pi_{k-1} F_{s} \map \cdots
\]
of abelian groups.  These sequences assemble to give the 
desired sequence.
\end{proof}

\begin{proposition}
\label{prop:pro-homotopy-cofibration}
Suppose that $j: A \map B$ is a cofibration of pro-spectra
such that $\pi_k j$ is an isomorphism of pro-groups for every $k$
and such that $j$ is an essentially levelwise $n$-equivalence
for some $n$.  Then $j$ is a $\pi_*$-acyclic cofibration.
\end{proposition}

\begin{proof}
We may assume that $j$ is a
levelwise $n$-cofibration because of Proposition 
\ref{prop:level-n-cofibration}.
We will show that $j$ is a levelwise $(n+1)$-cofibration.
By induction, this will imply that $j$ is a levelwise
$m$-cofibration for every $m$ and hence a levelwise $m$-equivalence
for every $m$ by Proposition \ref{prop:n-cofibration}.
This means that $j$ is a $\pi_*$-weak equivalence.

For any $s$, we have a map $j_s: A_s \map B_s$.  Factor $j_s$
into an $(n+1)$-cofibration $i_{s,n+1}: A_s \map Z_{s,n+1}$ 
followed by a co-$(n+1)$-fibration $p_{s,n+1}: Z_{s,n+1} \map B_s$
as in Definition \ref{defn:F}.

Let $C$ be the cofiber of $j$, which we may assume is constructed
levelwise.
From the long exact sequence of Proposition \ref{prop:cofiber-sequence},
we see that $\pi_{n} C$ is the trivial pro-group.
Therefore, we may choose $t \geq s$ such that the map 
$\pi_{n+1} C_{t} \map \pi_{n+1} C_{s}$ is zero.
Note also that the map $* \map C_t$ is an $n$-equivalence
because the map $A_t \map B_t$ is an $n$-cofibration.

Let $F$ be the fiber of $p_{s,n+1}$.  Note that the map $F \map *$
is a co-$(n+1)$-equivalence because $p_{s,n+1}$ is 
a co-$(n+1)$-fibration.

Consider the diagram
\[
\xymatrix{
A_t \ar[r]\ar[d] & A_s \ar[r]\ar[d] & Z_{s,n+1} \ar[d] \\
B_t \ar[r] & B_s \ar[r] & B_s.                           }
\]
The obstruction to lifting the right square is an element $\alpha$
of $[\Omega C_s, F]$, and the obstruction to lifting the outer
rectangle is the image of $\alpha$ under the map
$[\Omega C_s, F] \map [\Omega C_t, F]$ \cite[Cor.~8.4]{CDI}.  
We have chosen $t$ such that 
the map $\pi_n \Omega C_t \map \pi_n \Omega C_s$
is zero.  Also, note that $* \map \Omega C_t$ is 
an $(n-1)$-equivalence because its suspension
$* \map C_t$ is an $n$-equivalence.
Therefore, the conditions of Lemma \ref{lem:zero-map}
apply, and we conclude that the map 
$[\Omega C_s, F] \map [\Omega C_t, F]$ is zero.
Thus, the obstruction for lifting the outer square, which lies
in the image of this map, must be zero, and a lift $h$ exists.

Using this lift $h$, we get a diagram
\[
\xymatrix{
A_t \ar[r] \ar[d] & A_s \ar[d]\ar[r] & A_s \ar[d] \\
B_t \ar[r]_-h & Z_{s,n+1} \ar[r] & B_s.            }
\]
Now the conditions of Lemma \ref{lem:pro-isomorphism} are satisfied,
so the objects $Z_{s,n+1}$ assemble into a pro-spectrum that is
isomorphic to $B$.  The maps $A_s \map Z_{s,n+1}$ are thus
a level representation for $j$; this demonstrates that $j$
is an essentially levelwise $(n+1)$-equivalence.
\end{proof}

\begin{theorem}
\label{thm:pro-homotopy}
A map of pro-spectra $f$ is a $\pi_*$-weak equivalence
(see Definition \ref{defn:weak-equivalence}) if and only if
$\pi_{k} f$ is an
isomorphism of pro-abelian groups for every $k$ 
and $f$ is an essentially levelwise $n$-equivalence for some $n$.
\end{theorem}

The second condition in the above theorem feels unnatural.
It sounds plausible to construct a model structure on pro-spectra
in which the weak equivalences are just pro-homotopy group 
isomorphisms, but we have no idea how to do this.
One way to rationalize the existence of the second condition is that 
the cofibrations and fibrations are both plausible, and
this leaves no choice in what the weak equivalences are.

\begin{proof}
First suppose that $f$ is a $\pi_*$-weak equivalence,
so $f$ is an essentially levelwise $n$-equivalence
for every $n$.  For any $k$, choose a level representation for
$f$ that is a levelwise $(k+1)$-equivalence.  Then $\pi_k f$
is a levelwise isomorphism, so it is an isomorphism of pro-groups.   This
finishes one implication.

For the other implication, suppose that
$\pi_{k} f$ is an
isomorphism of pro-abelian groups for every $k$ 
and $f$ is an essentially levelwise $n$-equivalence for some $n$.
Factor $f$ into a
cofibration $i$ followed by a strict acyclic fibration $p$.
Therefore, $p$ is a strict weak equivalence.  This means
that $\pi_k p$ is an essentially levelwise isomorphism, so it is an
isomorphism of pro-groups.  We can conclude that $\pi_k i$ is an isomorphism
of pro-groups.

Also note that $i$ is an essentially levelwise 
$(n-1)$-equivalence by Lemma \ref{lem:2/3-n-equivalence}.
Therefore, Proposition \ref{prop:pro-homotopy-cofibration} applies,
and we can conclude that $i$ is a $\pi_*$-weak equivalence.
Since $p$ is also a $\pi_*$-weak equivalence,
the two-out-of-three axiom tells us that $f$ is
a $\pi_*$-weak equivalence.
\end{proof}

\section{Cohomology and the Whitehead Theorem} \label{sctn:Whitehead}

One of the primary motivations for the construction of our model structure
is the study of cohomology of pro-spectra.  
We now explore the relationship between 
$\pi_*$-weak equivalences and cohomology isomorphisms.  
We recall first the definition of cohomology for pro-spectra.

Let $HA$ be a fibrant Eilenberg-Mac Lane spectrum such that
$\pi_0 HA = A$, where $A$ is an abelian group.

\begin{definition} \label{defn:cohomology}
The \mdfn{$r$th cohomology $H^{r}(X;A)$ with coefficients in $A$}
of a pro-spectrum $X$ is the abelian group 
$[X, cHA]^{-r}_{\pros}$.
\end{definition}

This definition is precisely analogous to the definition of ordinary
cohomology for spectra.  

\begin{proposition} \label{prop:ordinary-cohomology}
If $X$ is any pro-spectrum, then
the cohomology group $H^{r} (X; A)$ is isomorphic to
$\colim_{s} H^{r}(X_{s}; A)$.
\end{proposition}

\begin{proof}
The spectrum $HA$ has no homotopy groups above dimension $0$, so
Proposition \ref{prop:into-constant} applies.
\end{proof}

The previous proposition shows that our definition of ordinary
cohomology in terms of Eilenberg-Mac Lane constant pro-spectra agrees
with the traditional notion of the cohomology of a pro-object.
Our viewpoint is that the straightforward colimit formula for
cohomology works because $HA$ is bounded above and 
therefore $cHA$ is $\pi_*$-fibrant.

We now work toward a Whitehead theorem for detecting $\pi_*$-weak
equivalences in terms of cohomology.

\begin{lemma}
\label{lem:lift-EM}
Let $i:A \map B$ be a cofibration that is an ordinary cohomology isomorphism
for all coefficients,
and let $X \map Y$
be a fibration of spectra whose fiber $F$ is an Eilenberg-Mac Lane
spectrum.  Then $i$ has the left lifting property
with respect to the constant map $q: cX \map cY$.
\end{lemma}

\begin{proof}
We may assume that $i$ is a levelwise cofibration.  Let $C$ be the 
cofiber of $i$, which we may assume is constructed levelwise.
The long exact sequence in cohomology for a cofiber sequence
indicates that the cohomology of $C$ is zero.
Let $k$ be the integer such that $\pi_k F$ is non-zero.

Consider a square
\[
\xymatrix{
A \ar[d]_i \ar[r] & cX \ar[d]^q \\
B \ar[r] & cY                }
\]
of pro-spectra.  This diagram is represented by a square
\[
\xymatrix{
A_s \ar[d]_{i_s} \ar[r] & X \ar[d]^q \\
B_s \ar[r] & Y                }
\]
of spectra.  The obstruction $\alpha$ to lifting this square is a weak
homotopy class belonging to $[\Omega C_s, F]$ \cite[Rem.~8.3]{CDI}.

Note that $[\Omega C_s, F]$ equals $H^{k+1}(C_s; \pi_k F)$
because $F$ is an Eilenberg-Mac Lane spectrum.  
Because $H^{k+1}(C; \pi_k F)$ is zero, there exists a $t$ such
that $\alpha$ pulls back to zero in $H^{k+1}(C_t;\pi_k F)$.  Therefore,
the obstruction to lifting the square
\[
\xymatrix{
A_t \ar[r] \ar[d]_{i_t} & A_s \ar[r] & X \ar[d]^q \\
B_t \ar[r] & B_s \ar[r] & Y                }
\]
vanishes, and a lift exists.  This lift represents the desired lift.
\end{proof}

\begin{theorem} \label{thm:Whitehead}
Let $f: X \map Y$ be an essentially levelwise $n$-equivalence for some $n$.
Then $f$ is a $\pi_*$-weak equivalence if and only if it is an
ordinary cohomology isomorphism for all coefficients.
\end{theorem}

\begin{proof}
One direction is easy; since cohomology is represented in the
$\pi_*$-homotopy category, $\pi_*$-weak equivalences are cohomology
isomorphisms.

For the other direction,
let $f: X \map Y$ be an essentially levelwise $n$-equivalence for some $n$
and an ordinary cohomology isomorphism for all coefficients.
Factor $f$ into a cofibration $i:X \map Z$ 
followed by a $\pi_*$-acyclic fibration $p:Z \map Y$.
Then $i$ is still an essentially levelwise $n$-equivalence for some $n$
by Lemma \ref{lem:2/3-n-equivalence}.  
Also, $p$ is a $\pi_*$-weak equivalence, so it
is a cohomology isomorphism since cohomology is defined to be representable.
This means that $i$ is a cohomology isomorphism, and we just have to
show that $i$ is a $\pi_*$-acyclic cofibration.

We may assume that $i$ is a levelwise $n$-cofibration by 
Proposition \ref{prop:level-n-cofibration}.
We use Lemma \ref{lem:detect-acyclic-cofibrations} to show
that $i$ is a $\pi_*$-acyclic cofibration.  Thus, we must find a lift
in the square
\[
\xymatrix{
X \ar[d]_i \ar[r] & cE \ar[d] \\
Z \ar[r] & cB                }
\]
of pro-spectra,
where $q:E \map B$ is a co-$m$-fibration for some $m$.
This diagram is represented by a square
\[
\xymatrix{
X_s \ar[d]_{i_s} \ar[r] & E \ar[d] \\
Z_s \ar[r] & B                }
\]
of spectra, and we have to find a lift after refining $s$.
According to Lemma \ref{lem:co-n-fibration-composition},
$q$ is a retract of a finite composition of maps that are
either co-$n$-fibrations or fibrations whose fiber is an Eilenberg-Mac Lane
spectrum.  Since we are trying to solve a lifting problem, we 
may assume that $q$ is a co-$n$-fibration or has an Eilenberg-Mac Lane
spectrum as its fiber.

If $q$ is a co-$n$-fibration, then a lift exists without refining $s$
at all since $i_s$ is an $n$-cofibration.
In the other case, Lemma \ref{lem:lift-EM} produces the lift.
\end{proof}

\section{Atiyah-Hirzebruch Spectral Sequence for Pro-Spectra} 
\label{sctn:AHSS}

We now consider generalized cohomology for pro-spectra.

\begin{definition} \label{defn:generalized-cohomology}
Let $E$ be any fixed pro-spectrum.
The \mdfn{$r$th $E$-cohomology $E^{r}(X)$}
of a pro-spectrum $X$ is the abelian group 
$[X, E]^{-r}_{\pros}$.
\end{definition}

This definition is precisely analogous to the definition of generalized
cohomology for spectra.  
In general, the calculation of $E^{*}X$ requires a fibrant replacement
of the pro-spectrum $E$.
When $E$ is constant, the Postnikov tower of $E$ is
one possible such fibrant replacement.

For example, if $E$ is the constant pro-spectrum $cKU$, then
$KU^r(X)$ is equal to $[X, P_* KU]_{\pros}$, where $P_* KU$ is
the Postnikov tower of $KU$.  Since $KU$ is not bounded above,
it is not true that $KU^r(X)$ is equal to $\colim_s KU^r(X_s)$;
the hypothesis of Proposition \ref{prop:into-constant} is not
satisfied.

We now develop an analogue of the Atiyah-Hirzebruch spectral sequence
\cite{AH}
for pro-spectra.  Here we are addressing the question of computing
$[X, Y]^{-r}_{\pros}$ for an arbitrary pro-spectrum $X$ and an arbitrary
pro-spectrum $Y$.  

We fix a pro-spectrum $Y$.
Let $A^q$ be the pro-spectrum $P_{-q} Y$, i.e., the levelwise
$(-q)$th Postnikov section of $Y$.
We choose the unusual indexing on $A^{q}$ in order to standardize the
indexing of our cohomological spectral sequence.
There is a diagram
\[
\cdots \map A^{q-1} \map A^{q} \map A^{q+1} \map \cdots,
\]
and we let $A$ be the inverse limit (in the category of pro-spectra)
of this tower.  As a cofiltered
diagram, $A$ is described by $(s,q) \mapsto P_{-q} Y_s$.
Note that a strict fibrant replacement for $A$ is a 
$\pi_*$-fibrant replacement for $Y$.

For any spectrum $Z$, let $C_q Z$ be the $q$th connected cover of $Z$,
i.e., the homotopy fiber of the map $Z \map P_q Z$.
Let $B^q$ be the pro-spectrum 
$C_{-q} Y$, i.e., the levelwise $(-q)$-connected cover of $Y$.
Again, there is a diagram
\[
\cdots \map B^{q-1} \map B^{q} \map B^{q+1} \map \cdots,
\]
and we let $B$ be the inverse limit (in the category of pro-spectra)
of this tower.  As a cofiltered
diagram, $B$ is described by $(s,q) \mapsto C_{-q} Y_s$.

We next show that $B$ is contractible in the $\pi_*$-model structure.

\begin{lemma} \label{lem:contractible}
The map $* \map B$ is a $\pi_*$-weak equivalence.
\end{lemma}

\begin{proof}
Fix an integer $n$.
The pro-spectrum $B = \lim_q B^{q}$ is isomorphic to the pro-spectrum
$\lim_{q < -n} B^q$.  Every object of $\lim_{q<-n} B^q$ 
is of the form $C_{-q} Y_s$ for some $s$ and some $q<-n$.  The map
$* \map C_{-q} Y_s$ is an $n$-equivalence because $-q > n$.
Thus, $* \map \lim_{q< -n} B^q$ is a levelwise $n$-equivalence,
so $* \map B$ is an essentially levelwise $n$-equivalence.  Since $n$
was arbitrary, this shows that $* \map B$ is a $\pi_*$-weak equivalence.
\end{proof}

Recall that for every $s$, $\Sigma^{-q} H\pi_{-q} Y_s$ 
is an Eilenberg-Mac Lane spectrum whose only non-zero homotopy group lies
in dimension $-q$ and is isomorphic to $\pi_{-q} Y_s$.
Let
$\Sigma^{-q} H\pi_{-q} Y$ be the obvious pro-spectrum 
constructed out of these Eilenberg-Mac Lane spectra.

\begin{lemma} \label{lem:homotopy-fiber}
For every $q$, the sequence
\[
B^q \map B^{q+1} \map \Sigma^{-q} H\pi_{-q} Y
\]
is a homotopy cofiber sequence of pro-spectra.
\end{lemma}

\begin{proof}
In order to compute the 
homotopy cofiber of any map, we should replace it by a levelwise
cofibration and then take the cofiber, i.e., the levelwise cofiber.
In other words, we just need to take the levelwise homotopy cofiber.

Recall that $B^q \map B^{q+1}$ is given levelwise by maps
$C_{-q} Y_s \map C_{-q-1} Y_s$.  
The homotopy cofiber of $C_{-q} Y_s \map C_{-q-1} Y_s$
is $\Sigma^{-q} H\pi_{-q} Y_s$.
\end{proof}

Let $X$ be any pro-spectrum.
Define \mdfn{$D_{2}^{p,q}$} to be $[X, B^{q}]^{p+q}_{\pros}$, and define
\mdfn{$E_{2}^{p,q}$} to be $[X, \Sigma^{-q}H\pi_{-q}Y ]^{p+q}_{\pros}$.  
The $\pi_*$-model structure is stable from Theorem \ref{thm:pro-stable},
so the homotopy cofiber sequence of
Lemma \ref{lem:homotopy-fiber} is also a homotopy fiber
sequence \cite[Thm.~7.1.11]{Ho}.
After applying the functor $[X, -]_{\pros}$, one obtains a
long exact sequence.  Therefore, we have an exact couple
\[
\xymatrix{
D_{2} \ar[rr]^{(-1,1)} & & D_{2} \ar[dl]^{(1,-1)} \\
& E_{2} \ar[ul]^{(1,0)}       }
\]
in which the labels indicate the degrees of the maps.
A careful inspection of degrees shows that this gives us a spectral
sequence beginning with the $E_{2}$-term.

Now we have a spectral sequence, but we must study its convergence.
We take the viewpoint of \cite{B}.

\begin{lemma} \label{lem:convergence1}
For all $n$ and all $X$,
the groups $\lim_{q} D^{q,n-q}_{2}$ and
$\lim_{q}^1 D^{q,n-q}_{2}$ vanish as $q \map \infty$.
\end{lemma}

\begin{proof}
Let $\hat{B}^q$ be the pro-spectrum 
described by $(s,p) \mapsto P_{-p} C_{-q} Y_s$.  The map
$B^q \map \hat{B}^q$ is a $\pi_*$-weak equivalence and $\hat{B}^q$
is levelwise bounded above.  By taking a levelwise fibrant replacement,
we may additionally assume that $\hat{B}^q$ is levelwise fibrant.

Let $\hat{B}$ be the inverse limit $\lim_q \hat{B}^q$ (computed in the
category of pro-spectra).  As a cofiltered diagram, $\hat{B}$ is described
by $(s,p,q) \mapsto P_{-p} C_{-q} Y_s$.
Again, $\hat{B}$ is levelwise bounded above, and the map $B \map \hat{B}$
is a $\pi_*$-weak equivalence.  Since each $\hat{B}^q$ is levelwise fibrant,
so is $\hat{B}$.

We already know that $[X,B]^{n}_{\pros}$ is zero for all $n$
because $B$ is contractible.
However, we will compute these homotopy classes 
another way by considering the components of the appropriate 
homotopically correct mapping space.

Take a cofibrant model for $\Sigma^{-n} X$.  
Corollary \ref{cor:fibrant2} says that the homotopically correct
mapping space for computing maps from $\Sigma^{-n} X$ to $B$
is weakly equivalent to 
$\holim_{s,p,q} \colim_t \Map(\Sigma^{-n} X_t, \hat{B}^q_{s,p})$.
Since homotopy limits commute, we can compute this as
\[
\holim_q \holim_{s,p} \colim_t \Map(\Sigma^{-n} X_t, \hat{B}^q_{s,p}).
\]
Now for a fixed $q$,
\[
M^q = \holim_{s,p} \colim_t \Map(\Sigma^{-n} X_t, \hat{B}^q_{s,p})
\]
is weakly equivalent again by Corollary \ref{cor:fibrant2}
to the homotopically correct mapping space for computing
maps from $\Sigma^{-n} X$ to $B^q$.

Now apply the short exact sequence of \cite[Thm.IX.3.1]{BK}
for computing the homotopy groups of the homotopy limit of the countable
tower
\[
\cdots \map M^{q-1} \map M^q \map M^{q+1} \map \cdots.
\]
We obtain the sequence
\[
0 \map \lim\nolimits^1_q \pi_1 M^q \map
\pi_0 \holim\nolimits_q M^q \map \lim\nolimits_q \pi_0 M^q \map 0.
\]
We already know that the middle term of the sequence is zero because
it is equal to $[X, B]^n_{\pros}$.
Therefore, the first and last terms are also
zero.  The first term is $\lim^1_q [X, B^q]^{n-1}_{\pros}$, and the
last term is $\lim_q [X,B^q]^{n}$.

Since $n$ is arbitrary, we have shown that as $q \map -\infty$, both
$\lim_q [X,B^q]^n$ and $\lim^1_q [X,B^q]^n$ are zero for all $n$.  Changing
indices gives that 
$\lim_q [X,B^{n-q}]^n$ and $\lim^1_q [X,B^{n-q}]^n$ are both zero
as $q \map \infty$.
\end{proof}

\begin{lemma} \label{lem:convergence3}
Suppose that the map $* \map X$ is an essentially 
levelwise $m$-equivalence for some $m$.
For all $n$, as $q \map - \infty$, $\colim_{q} D^{q,n-q}_{2}$ 
is isomorphic to $[X, Y]^{n}_{\pros}$.
\end{lemma}

\begin{proof}
We want to show that the natural map
\[
\colim_{q} [X, B^{n-q}]_{\pros}^{n} \map [X, Y]_{\pros}^{n}
\]
is an isomorphism.
The sequence $B^{n-q} \map Y \map A^{n-q}$ is a levelwise homotopy cofiber
sequence (given by $C_{q-n} Y_s \map Y_s \map P_{q-n} Y_s$ for each $s$), 
so it is a homotopy cofiber sequence of pro-spectra.  Therefore,
it is also a homotopy fiber sequence, so we get a long exact
sequence after applying $[X, -]_{\pros}$.  By the usual argument
with long exact sequences and the exactness of filtered colimits,
it suffices to show that $\colim_{q} [X, A^{n-q}]_{\pros}^{n}$ is 
zero for every $n$.

We may assume that the map $* \map X$ is a levelwise $m$-equivalence,
so the map $* \map \Sigma^{-n} X$ is a levelwise $(m-n)$-equivalence.

Choose $q' = m-1$, so $q'-n = m-n-1$.
Then the map $A^{n-q'} \map *$ is a levelwise co-$(m-n)$-equivalence
because each object of $A^{n-q'}$ is an $(m-n-1)$st Postnikov section.
Now the hypotheses of Corollary \ref{cor:pro-connected-Postnikov}
are satisfied, so $[\Sigma^{-n} X, A^{n-q'}]_{\pros}$ is zero.
Therefore, the colimit $\colim_q [X, A^{n-q}]_{\pros}^n$ is zero.
\end{proof}

\begin{lemma} \label{lem:convergence2}
Suppose that $X$ is an essentially bounded below pro-spectrum
and that $Y$ is a constant pro-spectrum.
For all $n$, as $q \map - \infty$, $\colim_{q} D^{q,n-q}_{2}$ 
is isomorphic to $[X, Y]^{n}_{\pros}$.
\end{lemma}

The hypothesis means that $X$ is isomorphic to a pro-spectrum $X'$
such that every $X'_{s}$ is bounded below, but the dimension at which
the homotopy groups of $X'_{s}$ vanish may depend on $s$.
For example, the pro-spectra $\R P^{\infty}_{-\infty}$ and
$\Cx P^{\infty}_{-\infty}$ \cite{L} are
essentially bounded below.

\begin{proof}
We want to show that the natural map
\[
\colim_{q} [X, B^{n-q}]^n_{\pros} \map [X, Y]^n_{\pros}
\]
is an isomorphism.
As in the proof of Lemma \ref{lem:convergence3},
it suffices to show that the group 
$\colim_{q} [\Sigma^{-n} X, A^{n-q}]_{\pros}$ is 
zero for every $n$. 

We may assume that $X$ is levelwise bounded below; then $\Sigma^{-n} X$
is also levelwise bounded below.  

By Proposition \ref{prop:into-constant}, we must show that
$\colim_{q,t} [\Sigma^{-n} X_{t}, A^{n-q}]$ is zero; here is where we use that
$Y$ and therefore $A^{n-q}$ is constant.
Fix an index $t$.
By the assumption on $X$,
there exists $m$ such that $* \map X_t$ is an $m$-equivalence.
Then the map $* \map \Sigma^{-n} X_t$ is an $(m-n)$-equivalence.

Choose $q' = m$, so $q' - n = m - n$.
Now the homotopy groups of $A^{n-q'}*$ 
vanish in dimensions greater than or equal to $m-n+1$ because
$A^{n-q'}$ is an $(m-n)$th Postnikov section.
The homotopy groups of $\Sigma^{-n} X_t$ vanish in dimensions less
than or equal to $m-n$.
These conditions on the homotopy groups guarantee that 
$[\Sigma^{-n} X_{t}, A^{n-q'}]$ is zero.
This shows that the colimit is zero.
\end{proof}

\begin{theorem} \label{thm:AHSS}
Let $X$ and $Y$ be any pro-spectra.
There is a spectral sequence
\[
E_{2}^{p,q} = H^{-p} (X; \pi_{-q} Y)
\]
converging to $[X, Y]_{\pros}^{p+q}$.
The differentials have degree $(r, -r + 1)$.
The spectral sequence is conditionally convergent if:
\begin{enumerate}[(1)]
\item
$X$ is essentially levelwise bounded below, or
\item
$Y$ is a constant pro-spectrum and $* \map X$ is an essentially levelwise
$n$-equivalence for some $n$.
\end{enumerate}
\end{theorem}

\begin{proof}
The conditional convergence comes from Lemmas \ref{lem:convergence1},
\ref{lem:convergence3}, and \ref{lem:convergence2}.
The identification of the $E_{2}$-term is given in 
Definition \ref{defn:cohomology}.
\end{proof}

\bibliographystyle{amsalpha}

\end{document}